\documentclass[reqno]{amsart}

\usepackage{amsmath, amssymb, amsthm, mathtools,xcolor,enumitem, parskip}
\usepackage{graphicx}
\usepackage{comment}
\usepackage{setspace}
\usepackage[top=1.5in, bottom=1.5in, left=1.1in, right=1.1in]{geometry}
\usepackage{hyperref}
\hypersetup{
  colorlinks=true,
  linkcolor=blue,
  citecolor=blue,
  urlcolor=blue,
  linktoc=all,
  bookmarksopen=true,
  bookmarksnumbered=true
}

\newtheorem{theorem}{Theorem}[section]
\newtheorem{lemma}[theorem]{Lemma}
\newtheorem{corollary}[theorem]{Corollary}
\newtheorem{proposition}{Proposition}
\newtheorem{conjecture}[theorem]{Conjecture}
\theoremstyle{definition}
\newtheorem{definition}{Definition}

\newtheorem{remark}{Remark}

\theoremstyle{definition}

\DeclareMathOperator{\Crit}{Crit}

\title[Relatively \'{E}tale Systems over Products of Curves]{The Dynamical Mordell--Lang Conjecture for Relatively \'{E}tale Systems over Products of Curves}

\author{Bhawesh Mishra}

\address{(Currently Unaffiliated) Previously at: 384 Dunn Hall, Department of Mathematical Sciences, The University of Memphis, Memphis, TN 38152}
\curraddr{Ramjanki Path, Biratnagar, Morang, Kosi, Nepal 56613}
\email{bhaweshmishra2024@gmail.com}

\begin{document}


\date{\today}

\begin{abstract}
We establish the dynamical Mordell--Lang conjecture over \(\mathbb C\) for endomorphisms that are relatively \'{e}tale over arbitrary endomorphisms of finite products of smooth projective curves. In particular, we establish the conjecture for arbitrary endomorphisms of product of smooth projective curves and relatively \'{e}tale skew products. The arithmetic input is our unconditional theorem for split endomorphisms over \(\overline{\mathbb Q}\), proved by reducing to a simultaneous Hasse principle for the local periods of critical points and analyzing local inertia in joint arboreal towers. We combine this arithmetic input with relative \'{e}taleness, simultaneous algebraic specialization preserving wandering base coordinates, prescribed-reduction \(p\)-adic embeddings, and \(p\)-adic interpolation to obtain the result in the complex case.
\end{abstract}

\maketitle
\tableofcontents

\part{Introduction and Summary of Our Methods}

\section{Introduction}
The dynamical Mordell--Lang conjecture embodies a remarkable convergence of arithmetic geometry, algebraic dynamics, Diophantine approximation, and combinatorics around the arithmetic structure of orbit--subvariety intersections. Given a self-map \(F:X\to X\) of varieties, a point \(x\in X\), and a closed subvariety \(V\subseteq X\), the dynamical Mordell--Lang conjecture predicts that \(
S_F(x,V)=\{n\ge0:F^n(x)\in V\}\) is a finite union of arithmetic progressions. The goal of this article is to confirm the conjecture over several broad classes of settings. We fix the notation that \(\mathbb{N}_{0} = \{0, 1, 2, 3, \ldots \}\). 

\begin{definition}
    Suppose that \(X\) and \(B\) are varieties, \(p:X\to B\), \(F:X\to X\), and \(\beta:B\to B\) are morphisms satisfying \(p\circ F=\beta\circ p\). The relative graph morphism of \(F\) over \(\beta\) is \(\delta_F:X\to B\times_{\beta,B,p}X\), given by \(x\mapsto(p(x),F(x))\), where the fiber product uses \(\beta\) on the first factor and \(p\) on the second. We call \(F\) \emph{relatively \'{e}tale over \(\beta\)} if \(\delta_F\) is \'{e}tale. 
\end{definition}

\begin{remark}
On the fiber over \(b\), \(F\) induces \(F_b:X_b\longrightarrow X_{\beta(b)},
\) and relative \'{e}taleness means that \(F\) introduces no ramification beyond whatever ramification is already present in \(\beta\). For example, if \(B\) is a point, relative \'{e}taleness is ordinary \'{e}taleness of \(F\). Similarly, if \(X = B\), \(p=\operatorname{id}\), and \(F=\beta\), the relative graph is an isomorphism.
\end{remark}
Our most general result establishes the conjecture for relatively-\'{e}tale systems over systems of product of curves as stated below. 

\begin{theorem}[Relatively \'{e}tale systems over products of curves]
\label{thm:relative-curve-base}
Let \(B=\prod_{i=1}^{s}C_i\) be a finite product of smooth projective connected curves over \(\mathbb C\), with the empty product understood as a point, and let \(\beta:B\to B\) be any endomorphism. Let \(X\) be a quasi-projective variety over \(\mathbb C\), let \(p:X\to B\) be smooth, and let \(F:X\to X\) satisfy \(p\circ F=\beta\circ p\). If \(F\) is relatively \'{e}tale over \(\beta\), then, for every \(x\in X(\mathbb C)\) and every closed subvariety \(V\subseteq X\), the set \(\{n\in\mathbb N_0:F^n(x)\in V\}\) is a finite union of arithmetic progressions.
\end{theorem}

Theorem \ref{thm:relative-curve-base} confirms validity of the dynamical Mordell--Lang conjecture for several broad classes of systems, some of which now follow. Taking
\(X=B\), \(p\) to be the identity, and \(F=\beta\) in Theorem
\ref{thm:relative-curve-base} gives Theorem
\ref{thm:all-curve-products} below.

\begin{theorem}[Relatively \'{e}tale systems over split rational bases]\label{thm:relative-complex}
Let \(g\geq0\), let \(B=(\mathbb P^1_{\mathbb C})^g\), and let \(\beta=(\varphi_1,\ldots,\varphi_g):B\to B\), where every \(\varphi_i\) has degree at least \(2\); for \(g=0\), interpret \(B\) as a point. Let \(X\) be a quasi-projective variety over \(\mathbb C\), let \(p:X\to B\) be smooth, and let \(F:X\to X\) satisfy \(p\circ F=\beta\circ p\). If \(F\) is relatively \'{e}tale over \(\beta\), then, for every \(x\in X(\mathbb C)\) and every closed subvariety \(V\subseteq X\), the set \(\{n\in\mathbb N_0:F^n(x)\in V\}\) is a finite union of arithmetic progressions.
\end{theorem}

The total map in Theorem \ref{thm:relative-complex} need not be \'{e}tale: all ramification forced by the split base is allowed, while the remaining directions are controlled by the relative graph. The theorem imposes no product decomposition on \(X\) or \(F\), and the dynamics in the fibers of \(p\) may vary with the base point. We will also establish the conjecture for arbitrary endomorphisms of arbitrary product of smooth projective connected curves.  

\begin{theorem}[Endomorphisms of products of curves]
\label{thm:all-curve-products}
Let \(X=\prod_{i=1}^{s}C_i\) be a finite product of smooth projective connected curves over \(\mathbb C\), with the empty product understood as a point, and let \(\beta:X \to X\) be any endomorphism. Then, for every \(x \in X(\mathbb{C})\) and every closed subvariety \(V\subseteq X\), the set \(\{n\in\mathbb{N} : \beta^{n}(x) \in V \}\) is a finite union of arithmetic progressions.
\end{theorem}

Recently, Zhang independently announced dynamical Mordell--Lang for split endomorphisms of products of curves in \cite{ZhangSplit2026}, which overlaps with the split-case of Theorem \ref{thm:all-curve-products}. Another immediate consequence of our methods is the confirmation of the dynamical Mordell--Lang conjecture over relatively \'{e}tale skew products as below. 

\begin{theorem}[Relatively \'{e}tale skew products]
\label{thm:skew}
Let \(B\) be a finite product of smooth projective connected curves over \(\mathbb C\), let \(\beta:B\to B\) be any endomorphism, let \(Y\) be smooth and quasi-projective over \(\mathbb C\), and let \(\psi:B\times Y\to Y\) be a morphism. Define \(F:B\times Y\to B\times Y\) by \(F(b,y)=(\beta(b),\psi(b,y))\). If \((b,y)\mapsto(b,\psi(b,y))\) is \'{e}tale, then the dynamical Mordell--Lang conjecture holds for \(F\).
\end{theorem}

Another short and sweet corollary of our results is as follows:

\begin{corollary}[Products of symmetric powers]
\label{cor:symmetric-powers}
For \(1\leq j\leq t\), let \(B_j\) be a finite product of smooth
projective connected curves, let \(p_j:X_j\to B_j\) be smooth with
\(X_j\) quasi-projective, and let \(\Theta_j:X_j\to X_j\) be an
endomorphism. Suppose that there are \(u_j\geq1\) and
\(\beta_j:B_j\to B_j\) such that
\[
p_j\circ\Theta_j^{u_j}=\beta_j\circ p_j
\]
and the relative graph
\[
X_j\longrightarrow
B_j\times_{\beta_j,B_j,p_j}X_j,
\qquad
x\longmapsto
\bigl(p_j(x),\Theta_j^{u_j}(x)\bigr)
\]
is \'etale. For arbitrary \(n_j\geq1\), the endomorphism
\[
\prod_{j=1}^t\operatorname{Sym}^{n_j}(\Theta_j)
\quad\text{of}\quad
\prod_{j=1}^t\operatorname{Sym}^{n_j}(X_j)
\]
satisfies dynamical Mordell--Lang. In particular, if each \(X_j\)
is a finite product of smooth projective connected curves and each
\(\Theta_j:X_j\to X_j\) is an arbitrary endomorphism, then this
conclusion holds for every choice of \(n_1,\ldots,n_t\geq1\).
\end{corollary}

We impose no relation among the degrees of the nonlinear rational factors, no genericity condition in a parameter space, no postcriticality, and no restriction on the dimension or position of \(V\) or the starting point \(x\) in our results. Certain weaker results regarding the dynamical Mordell--Lang conjecture permit an additional exceptional set of density zero, and gap principles show that a nonperiodic return set must be extremely sparse. Our results remove that possible infinite remainder. Every infinite set of return times contains an entire arithmetic-progression tail; equivalently, after passing to a residue class and an iterate, the relevant tail of the orbit is contained in \(V\). This is the return-time rigidity predicted by the conjecture itself.

\subsection{Motivation}

Mordell's conjecture, proved by Faltings in \cite{Faltings83}, asserts that a curve of genus greater than one over a number field has only finitely many rational points. Lang placed this finiteness statement in a broader intersection-theoretic setting by asking how a subvariety of a commutative algebraic group can meet a finitely generated subgroup. Faltings treated the case of abelian varieties in \cite{Faltings91,Faltings1994}; the semiabelian formulation was completed through the work of Vojta \cite{Vojta96} and McQuillan \cite{McQuillan}. In the form most closely related to dynamics, the resulting theorem is the following.

\begin{theorem}[Originally Mordell--Lang Conjecture]\label{thm:Mordell--Lang}
Let \(X\) be a semiabelian variety defined over complex numbers, let \(V \subseteq X\) be a subvariety and let \(\Gamma \subseteq X(\mathbb{C})\) be a finitely generated subgroup. Then, \(\Gamma \cap V(\mathbb{C})\) is a finite union of cosets of subgroups of \(\Gamma\).    
\end{theorem}

Theorem \ref{thm:Mordell--Lang} expresses one of the guiding principles of arithmetic geometry: an (irreducible) subvariety containing a Zariski-dense collection of special points must itself be special. In Theorem \ref{thm:Mordell--Lang}, special points lie in the finitely generated group \(\Gamma\), while the special irreducible subvarieties are translates of algebraic subgroups. If we restrict attention to the cyclic semigroup generated by one element, the set of exponents for which the orbit lies on \(V\) is a finite union of congruence classes. Dynamical Mordell--Lang asks for the same arithmetic regularity after a cyclic group orbit is replaced by the forward orbit of an arbitrary algebraic self-map. The modern formulation emerged from work of Ghioca and Tucker, notably \cite{Ghioca2009}; we use the following standard form.

\begin{conjecture}[\textbf{Dynamical Mordell--Lang Conjecture}]
Let \(\Phi\) be an endomorphism of a quasi-projective variety  \(X\) defined over \(\mathbb{C}\), \(V \subset X\) be a closed subvariety and \(x \in X(\mathbb{C})\). Then, the set of positive integers \(n\) such that \(\Phi^{n}(x) \in V(\mathbb{C})\) is a union of finitely many arithmetic progressions\(\{an+b\}_{n\in\mathbb{N}}\) with $a, b\geq 0$.
\label{con:DMLC}
\end{conjecture}

We allow a progression with common difference zero, so finite exceptional return times are absorbed into the same formulation. In this language, Conjecture \ref{con:DMLC} says that the indicator of the return set is eventually periodic. This is much stronger than a mere finiteness or density assertion: it excludes every infinite but arithmetically irregular pattern of intersections.

The following equivalent reformulation of the dynamical Mordell--Lang conjecture makes the philosophy of special points and special (irreducible) subvarieties explicit; see \cite[pp. 50-51]{DMLBook} for the equivalence.

\begin{conjecture}
\label{con:equivalentDMLC}
    Let \(X\) be a quasi-projective variety over \(\mathbb{C}\), \(\Phi\) be any endomorphism of \(X\), \(\alpha\in X(\mathbb{C})\) and \(V \subseteq X\) be an irreducible subvariety of positive dimension. If the intersection \(\mathcal{O}_{\Phi}(\alpha) \cap V(\mathbb{C})\) is Zariski dense in \(V\), then \(V\) is \(\Phi\)-periodic. 
\end{conjecture}

One direction is visible directly: an arithmetic progression of return times places a tail of an orbit inside a subvariety invariant under an iterate. Conversely, one analyzes the Zariski closures of the returns and applies Noetherian induction to the lower-dimensional pieces. This viewpoint situates the conjecture beside Manin--Mumford, Bogomolov, and Pink--Zilber, as well as their dynamical analogues: in each setting, an abundance of distinguished points forces geometric rigidity. These parallels are conceptual rather than automatic implications, but they explain why orbit--subvariety intersections have become a central organizing problem in arithmetic dynamics; see \cite{GhiocaTuckerZhang2011,GhiocaNguyenYe2019,BenedettoEtAl2019Survey}.

There is a second, complementary lineage. Let \(K\) be a field of characteristic zero, let \(r\geq1\), and let \((a_n)_{n\geq0}\) satisfy the fixed linear recurrence \(a_{n+r}=\sum_{j=0}^{r-1}c_ja_{n+j}\) for \(n\geq0\), where \(c_0,\ldots,c_{r-1}\in K\). After reducing the order and discarding finitely many initial terms when necessary, one may assume \(c_0\neq0\). The Skolem--Mahler--Lech theorem asserts that \(Z(a):=\{n\in\mathbb N_0:a_n=0\}\) is the union of a finite set and finitely many arithmetic progressions. Skolem proved the result for recurrences over \(\mathbb Q\) \cite{Skolem35}, Mahler extended it to recurrences over \(\overline{\mathbb Q}\) \cite{Mah35}, and Lech proved it over arbitrary fields of characteristic zero \cite{Lech53}. The connection with dynamics is given by the companion map \(f:\mathbb A_K^r\longrightarrow\mathbb A_K^r\), defined by \(f(y_0,\ldots,y_{r-1})=(y_1,\ldots,y_{r-1},\sum_{j=0}^{r-1}c_jy_j)\). For \(\alpha=(a_0,\ldots,a_{r-1})\), one has \(f^n(\alpha)=(a_n,a_{n+1},\ldots,a_{n+r-1})\). Consequently, if \(H\subseteq\mathbb A_K^r\) is the hyperplane \(y_0=0\), then \(a_n=0\) if and only if \(f^n(\alpha)\in H\). Thus the Skolem--Mahler--Lech theorem is the linear prototype of the dynamical Mordell--Lang conjecture. Bell's affine generalization developed this connection into an important dynamical tool \cite{Bell2006Affine,Bell2008Corrigendum}.

\subsection{Our Work in Context}\label{sec:previous-work}

The dynamical Mordell--Lang conjecture has attracted sustained attention because its simple statement conceals several different sources of difficulty: ramification of the map, failure of local analytic interpolation, the geometry of periodic subvarieties, and the possibility that infinitely many isolated returns might survive without forming a progression. The split case has long constituted a distinct central testing ground for
the dynamical Mordell--Lang conjecture; see, for example,
\cite[Ch. 7]{DMLBook}. One foundational approach is \(p\)-adic interpolation. Bell proved a generalized Skolem--Mahler--Lech theorem for automorphisms of affine varieties \cite{Bell2006Affine,Bell2008Corrigendum}; Bell, Ghioca, and Tucker then established dynamical Mordell--Lang for \'{e}tale endomorphisms of quasiprojective complex varieties \cite{Bell2010}. These results are strikingly general in the ambient variety and the target subvariety.  By Riemann--Hurwitz, every rational self-map of \(\mathbb{P}^{1}\) of degree at least two is ramified, and so a product of such maps is not \'{e}tale.

Other exact results proceed by methods adapted to a particular geometry. Fakhruddin proved strong orbit-density statements for generic endomorphisms of projective space \cite{Fakh14}. Xie first treated birational polynomial self-maps of \(\mathbb{A}^{2}\) \cite{Xiebirational} and later proved the full conjecture for all polynomial endomorphisms of the affine plane \cite{Xie17}; his valuation-theoretic work in \cite{Xie15} is an important ingredient in this circle of ideas. Ghioca and Xie established the conjecture for skew-linear self-maps \cite{GhiocaXieSkewLinear}. These theorems demonstrate that exact dynamical Mordell--Lang can be obtained well beyond the \'{e}tale setting, but they also show how strongly the known proofs depend on the special geometry of the system. Our result as in Theorem \ref{thm:relative-curve-base}, and Theorems \ref{thm:all-curve-products} and \ref{thm:skew}, extend this line of reasoning to a much broader class of previously unknown instances where the Conjecture \ref{con:DMLC} holds. 

A second line of work isolates what can be proved without resolving the exact return set. Bell, Ghioca, and Tucker showed that a continuous self-map of a Noetherian space has return sets equal to a finite union of arithmetic progressions together with a set of Banach density zero \cite{BGT15}; height-gap refinements and related weak forms appear in \cite{BellHuSatWeakDML}. For coordinatewise rational dynamics, Benedetto, Ghioca, Kurlberg, and Tucker established a gap principle showing that, if the return set contains no infinite arithmetic progression, then its elements grow extraordinarily rapidly \cite{BGKTGapPrinciple}. In this context, all our results close this gap and carry the exact conclusion as posited in the conjecture. 

The split polynomial cases supplied many of the first higher-dimensional testing grounds. Ghioca, Tucker, and Zieve related infinite intersections of polynomial orbits to common iterates and proved important dynamical Mordell--Lang cases for linear targets \cite{GTZ08,GTZ12}. Medvedev and Scanlon classified invariant subvarieties for coordinatewise polynomial systems, giving a structural account of the functional relations that can support periodic geometry \cite{MedvedevScanlon2014}. Their theorem classifies the special subvarieties; by itself, however, such a classification does not rule out an infinite sparse collection of nonperiodic returns. In a complementary direction, Ghioca and Nguyen established dynamical Mordell--Lang for curve targets under split polynomial maps in arbitrary dimension \cite{GhiocaNguyenBounds}. Thus the polynomial theory reached arbitrary numbers of factors when the target was a curve, and arbitrary targets in dimension two through Xie's theorem, but not arbitrary-dimensional targets for an arbitrary number of factors. 

For rational maps, even the diagonal case required new input. Benedetto, Ghioca, Kurlberg, and Tucker treated the diagonal action of a rational map on \((\mathbb P^1)^g\), with algebraic coefficients and curve targets, under the hypothesis that the map has no periodic superattracting points other than exceptional points; they also obtained further polynomial cases \cite{BGKT12}. More recent structural results clarify the rigidity behind exceptional orbit intersections. Pakovich classified broad families of invariant curves and proved common-iterate phenomena for tame rational functions \cite{Pakovich2023Curves,Pakovich2023Tame}, while Ghioca, Nguyen, and Ye proved dynamical Manin--Mumford and Bogomolov results for important classes of split rational maps \cite{GhiocaNguyenYe2019}. These results are closely related to the geometry predicted by Conjecture \ref{con:equivalentDMLC}, but they are not return-time theorems: identifying periodic subvarieties does not alone eliminate additional sparse intersections.

Genericity and arithmetic information about critical orbits have led to another substantial body of work. Juul used arboreal Galois information on backward critical orbits to obtain the split conclusion for typical polynomial maps in a precise parameter-space sense \cite{Juul}. These theorems cover large and natural families, yet their hypotheses exclude special algebraic relations among critical orbits.

The most recent results sharpen the boundary still further. Bell and Zhong proved \(p\)-adic interpolation for an orbit of an arbitrary rational map of \(\mathbb{P}^{1}\) and deduced dynamical Mordell--Lang when that rational factor is paired with an \'{e}tale system \cite{BellZhong2023}. Xie and Yang developed a height argument that, in characteristic zero, yields the full conjecture for split self-maps of products of curves when the nonlinear factor degrees are pairwise distinct \cite{XieYangHeight2025}. Xie, Yang, and Zheng proved the conjecture over \(\overline{\mathbb{Q}}\) for a product of an affine curve and a projective curve, a setting that includes a polynomial factor paired with a rational factor \cite{Xie26}. These advances show that neither rationality of an individual factor nor higher dimension alone is the decisive obstruction. What remained outside their scope was the unrestricted synchronization of repeated-degree, genuinely ramified rational factors---already for two factors---and, in higher products, arbitrary target dimension.

In addition to covering these previously separate regimes, our results also extend several of these lines of inquiry. All of our results permit every degree pattern among the factors of degree at least two, including complete repetition; every ramification and critical-orbit behavior; every positive number of factors; and every subvariety \(V\). In particular, they convert the gap principle from a sparsity dichotomy into the exact arithmetic-progression conclusion and extends beyond the curve targets that governed much of the earlier rational theory. Theorem \ref{thm:relative-complex} extends this arithmetic mechanism to arbitrary smooth systems relatively \'{e}tale over the product of curves, without requiring the total space or its endomorphism to split, while Theorem \ref{thm:relative-curve-base} allows arbitrary endomorphisms of products of curves as bases. For broader accounts of the conjecture and its place in arithmetic dynamics, we refer the readers to \cite{DMLBook,BenedettoEtAl2019Survey,XieSurvey2026}.
There also has been substantial interesting work towards establishing the dynamical Mordell--Lang conjecture in positive characteristic and interested readers can see \cite{PositiveCGSZ, PositiveSheXie, GhiocaPositive, CGLNpadic, LeeNamConverse} and the references therein.

\section{A Description of Our Approach}\label{sec:methods-overview}

Our approach can be summarized as having four steps. We first prove the conjecture for split rational maps over \(\overline{\mathbb Q}\). We then extend this arithmetic result to relatively \'{e}tale systems over split bases. Next, we pass from algebraic coefficients to arbitrary complex coefficients. Finally, we analyze stable images of endomorphisms of products of curves and reduce the general theorem to the split relative case.

\subsection{The arithmetic split case}

Our main arithmetic input is the following theorem.

\begin{theorem}[Arithmetic split theorem]\label{thm:split-arithmetic}
Let \(g\geq1\), and let \(\Phi=(\varphi_1,\ldots,\varphi_g)\) act coordinatewise on \((\mathbb P^1_{\overline{\mathbb Q}})^g\), where each \(\varphi_i:\mathbb P^1\to\mathbb P^1\) is defined over \(\overline{\mathbb Q}\) and has degree at least \(2\). Then, for every \(x\in(\mathbb P^1)^g(\overline{\mathbb Q})\) and every closed subvariety \(V\subseteq(\mathbb P^1_{\overline{\mathbb Q}})^g\), the set \(\{n\geq0:\Phi^n(x)\in V\}\) is a finite union of arithmetic progressions.
\end{theorem}

All the data in Theorem \ref{thm:split-arithmetic} are defined over a common number field, so it suffices to prove the theorem over number fields. Our \(p\)-adic input is the following result of Benedetto, Ghioca, Kurlberg, and Tucker \cite[Theorem~3.4]{BGKT12}.

\begin{theorem}\label{thm:padicattracting}
Let \(V\subseteq(\mathbb P^1)^g\) be a closed subvariety defined over \(\mathbb C_p\), and let \(\Phi=(\varphi_1,\ldots,\varphi_g)\), where every \(\varphi_i\in\mathbb C_p(t)\) has good reduction. Let \(\alpha=(x_1,\ldots,x_g)\in(\mathbb P^1(\mathbb C_p))^g\), and suppose that, for every \(i\), the orbit of \(x_i\) does not meet the residue class of any attracting \(\varphi_i\)-periodic point. Then \(V(\mathbb C_p)\cap\mathcal O_\Phi(\alpha)\) is a finite union of suborbits of the form \(\{\Phi^{kn+\ell}(\alpha):n\geq0\}\), where \(k,\ell\geq0\).
\end{theorem}

To apply Theorem \ref{thm:padicattracting}, we prove the following simultaneous nonperiodicity result for critical points.

\begin{theorem}[Simultaneous Hasse principle for critical points]\label{thm:simultaneousHasse}
Let \(K\) be a number field, let \(\varphi_1,\ldots,\varphi_g:\mathbb P^1_K\to\mathbb P^1_K\) have degree at least \(2\), and, for every \(i\), let \(\mathcal C_i=\{c_{ij}\}_{j=1}^{n_i}\subseteq\mathbb P^1(\overline K)\) be the set of nonperiodic critical points of \(\varphi_i\). Then there is a positive-density set of primes \(\mathfrak p\) of \(K\) such that, for every \(i,j\) and every prime \(\mathfrak q\) of \(K(c_{ij})\) above \(\mathfrak p\), the reduction of \(c_{ij}\) at \(\mathfrak q\) is not periodic under the reduction of \(\varphi_i\).
\end{theorem}

Section \ref{sec:inertia-orbits} proves that, above a periodic critical reduction, the number of relevant inertia orbits remains bounded while their lengths acquire arbitrarily large prime-power divisibility. Section \ref{sec:Group-combinatorics} uses this estimate to construct a derangement in a suitable Frobenius coset of a joint preimage field. Chebotarev's theorem then gives Theorem \ref{thm:simultaneousHasse} in Section \ref{sec:proof-Hasse-principle}.

In Section \ref{sec:Hasse-to-DML}, Theorem \ref{thm:simultaneousHasse} is strengthened to the critical-residue avoidance result of Theorem \ref{thm:intermediate}. It provides an integer \(M\) and a positive-density set of primes such that, after time \(M\), the chosen wandering orbits avoid the reductions of all critical points. Such an orbit cannot enter the residue class of an attracting periodic point. Theorem \ref{thm:padicattracting} therefore applies and proves Theorem \ref{thm:split-arithmetic}.

\subsection{Relative systems and the passage to \texorpdfstring{\(\mathbb C\)}{C}}

We next extend the arithmetic theorem to systems whose total space and total map need not split.

\begin{theorem}[Relatively \'{e}tale systems over split bases: the algebraic case]\label{thm:relative-algebraic}
Let \(K\) be a number field, let \(g\geq0\), let \(B=(\mathbb P^1_K)^g\), and let \(\beta=(\varphi_1,\ldots,\varphi_g):B\to B\), where every \(\varphi_i\) has degree at least \(2\); when \(g=0\), interpret \(B\) as a point. Let \(X\) be a quasi-projective \(K\)-variety, let \(p:X\to B\) be smooth, and let \(F:X\to X\) satisfy \(p\circ F=\beta\circ p\). If the relative graph \(\delta_F:X\to B\times_{\beta,B,p}X\) is \'{e}tale, then the dynamical Mordell--Lang conjecture holds for \(F_{\overline K}\).
\end{theorem}

The proof is given in Section \ref{sec:relative-systems}. After removing preperiodic base coordinates, Theorem \ref{thm:intermediate} supplies a prime at which the reduced base orbit avoids the ramification locus of \(\beta\). Hence the reduced base map is \'{e}tale along that orbit. Lemma \ref{lem:relative-etale-formal} transfers this orbitwise \'{e}taleness to the total map, and Lemma \ref{lem:orbitwise-interpolation} gives the required arithmetic-progression description of the return set.

Section \ref{sec:complex-transfer} passes from algebraic coefficients to arbitrary complex coefficients. Lemma \ref{lem:simultaneous-specialization} specializes the wandering base coordinates to algebraic wandering points and is used only to choose a favorable arithmetic reduction. Lemma \ref{lem:prescribed-reduction} then embeds the original coefficient field injectively into \(\mathbb Q_\ell\) while preserving that reduction. Relative \'{e}taleness and Lemma \ref{lem:orbitwise-interpolation} apply to the resulting \(p\)-adic system, while injectivity ensures that its return set is exactly the return set of the original complex system. This proves Theorem \ref{thm:relative-complex}.

It remains to reduce arbitrary endomorphisms of products of curves to this split situation. Section \ref{sec:stable-structure} proves the following structural theorem.

\begin{theorem}[Equivariant cover of the stable image]\label{thm:mixed-core}
Let \(B=\prod_{i=1}^gC_i\) be a product of smooth projective connected complex curves, and let \(F:B\to B\) be any endomorphism. There exist integers \(q\geq1\) and \(N\geq0\), a smooth projective variety \(Y\), a surjective \'{e}tale endomorphism \(\eta:Y\to Y\), morphisms \(\phi_1,\ldots,\phi_r:\mathbb P^1\to\mathbb P^1\) of degree at least \(2\), and a surjective morphism \(\lambda:(\mathbb P^1)^r\times Y\to B_0:=(F^q)^N(B)\) such that \(F^q(B_0)=B_0\) and \(F^q\circ\lambda=\lambda\circ(\phi_1,\ldots,\phi_r,\eta)\). No assertion is made that \(B_0\) is itself a product.
\end{theorem}

Theorem \ref{thm:mixed-core} follows from the fact that every nonconstant map from a product of curves to a curve of genus at least \(2\), and every nonconstant map from a product of projective lines to \(\mathbb P^1\), depends on only one factor. After passing to an iterate, the noncyclic dependencies disappear. The surviving nonlinear rational factors act coordinatewise, while the stable positive-genus part is covered by a surjective \'{e}tale system.

In Section \ref{sec:relative-curves-complex}, we pull the original relative system back through the cover supplied by Theorem \ref{thm:mixed-core}. Lemma \ref{lem:forget-etale-factor} absorbs the \'{e}tale factor \(Y\) into the relative directions, so Theorem \ref{thm:relative-complex} applies. Surjectivity transfers the conclusion to the stable image, and Lemma \ref{lem:iteration} passes from the chosen iterate and orbit tail to the original map. This proves Theorem \ref{thm:relative-curve-base}, and hence Theorems \ref{thm:all-curve-products} and \ref{thm:skew}.

We want to note an inherent limitation of our approach. Our methods rely essentially on criticality: the nontrivial local degrees at critical points produce the growing inertia orbits used in the Frobenius-coset argument. They currently do not yield an analogous simultaneous nonperiodicity theorem for arbitrary marked points. For the relative theorem, however, critical-point avoidance is precisely what is needed to make the reduced base map, and therefore the reduced total map, \'{e}tale along the chosen orbit. A random-model heuristic and numerical evidence appearing in \cite{BGHKST13} show that, for arbitrary points in dimension \(\geq 5\), this approach may not work. For the relative theorem, however, this is exactly the needed input: critical avoidance makes the split base map \'{e}tale along the reduced orbit, and relative \'{e}taleness transfers that orbitwise property to the total map. 

\subsection{Notation and terminology}
Given a number field \(K\) with algebraic closure \(\overline{K}\), \(\mathfrak{o}_{K}\) will denote the corresponding ring of algebraic integers. We will call a nonzero prime \(\mathfrak{p}\) of \(\mathfrak{o}_{K}\) simply a \textit{prime of \(K\)}, and we denote the corresponding residue field \(k_{\mathfrak{p}}:=\mathfrak{o}_{K}/\mathfrak{p}\). Upon fixing an isomorphism \(\pi\) from \({\mathbb{P}}^{1}_{K}\) to the generic fibre of \({\mathbb{P}}^{1}_{\mathfrak{o}_{K}}\), for each prime \(\mathfrak{p}\) of \(K\), and for each \(x\in{\mathbb{P}}^{1}(K)\), we denote by \(r_{\mathfrak{p}}(x)\) the intersection of the Zariski closure of \(\pi(x)\) with the fibre above \(\mathfrak{p}\) of \({\mathbb{P}}^{1}_{\mathfrak{o}_{K}}\). The resulting map \(r_{\mathfrak{p}}:{\mathbb{P}}^{1}(K)\longrightarrow{\mathbb{P}}^{1}(k_{ \mathfrak{p}})\) is the \textit{reduction map} at \(\mathfrak{p}\). We say that \(\alpha\in{\mathbb{P}}^{1}(K)\) is \textit{congruent} to \(\beta\in{\mathbb{P}^{1}}(K)\) modulo \(\mathfrak{p}\), and we write \(\alpha\equiv\beta\pmod{\mathfrak{p}}\), if \(r_{\mathfrak{p}}(\alpha)=r_{\mathfrak{p}}(\beta)\). For an \(\alpha\), we may also use \(\overline{\alpha}\) for its reduction modulo \(\mathfrak{p}\) when \(\mathfrak{p}\) is clear from context. We also use the following notion of good reduction of \(\varphi\), first introduced by Morton and Silverman in \cite{MS94}.

\begin{definition}\label{def:good-red}
Let \(K\) be a number field, let \(\mathfrak{p}\) be a prime of \(K\), and let \(\mathfrak{o}_{\mathfrak{p}}\subseteq K\) be the corresponding local ring of integers. Let \(\varphi:{\mathbb{P}}^{1}\longrightarrow{\mathbb{P}}^{1}\) be a morphism over \(K\), given by \(\varphi([X,Y])=[F(X,Y):G(X,Y)]\), where \(F,G\in\mathfrak{o}_{\mathfrak{p}}[X,Y]\) are relatively prime homogeneous polynomials of the same degree such that at least one coefficient of \(F\) or \(G\) is a unit in \(\mathfrak{o}_{\mathfrak{p}}\). Let \(\varphi_{\mathfrak{p}}:=[F_{\mathfrak{p}},G_{\mathfrak{p}}]\), where \(F_{\mathfrak{p}},G_{\mathfrak{p}}\in k_{\mathfrak{p}}[X,Y]\) are the reductions of \(F\) and \(G\) modulo \(\mathfrak{p}\). We say that \(\varphi\) has good reduction at \(\mathfrak{p}\) if \(\varphi_{\mathfrak p}:\mathbb P^1_{k_{\mathfrak p}}\to\mathbb P^1_{k_{\mathfrak p}}\) is a morphism of the same degree as \(\varphi\).
\end{definition}

We denote the reduced morphism in Definition \ref{def:good-red} by \(\overline\varphi\) when \(\mathfrak p\) is clear. If \(\varphi\in K[t]\) is a polynomial, there is an elementary criterion for good reduction: \(\varphi\) has good reduction at \(\mathfrak{p}\) if and only if all coefficients of \(\varphi\) are \(\mathfrak{p}\)-adic integers, and its leading coefficient is a \(\mathfrak{p}\)-adic unit. 

Throughout Sections \ref{sec:inertia-orbits}--\ref{sec:relative-systems}, \(K\) denotes a number field; Sections \ref{sec:complex-transfer}--\ref{sec:consequences} concern systems over \(\mathbb C\). Given a rational map \(\varphi:\mathbb{P}^{1}_{K}\longrightarrow\mathbb{P}^{1}_{K}\) of degree at least \(2\) and integer \(m\geq 1\), we write \(\varphi^{m}\) for the \(m\)-th iterate of \(\varphi\) under composition. The \textit{forward orbit} of a point \(\alpha\in\mathbb{P}^{1}(K)\) under \(\varphi\) is \(\mathcal{O}^{+}(\alpha)=\{\varphi^{m}(\alpha):m\geq 0\}\). We say that \(\alpha\) is \textit{\(\varphi\)-periodic} if \(\varphi^{m}(\alpha)=\alpha\) for some \(m\geq 1\), and \textit{\(\varphi\)-preperiodic} if its forward orbit is finite. A point that is not preperiodic is called \emph{wandering}.

Given a number field \(E\), a rational map \(\varphi : \mathbb{P}^{1}_{E} \longrightarrow \mathbb{P}^{1}_{E}\) of degree \(d \geq 2\), \(\alpha\in\mathbb{P}^{1}(E)\), and \(n \geq 1\), define \(E_{n}(\varphi, \alpha) := E\big(\varphi^{-n}(\alpha)\big)\). Thus \(E_n(\varphi,\alpha)\) is obtained by adjoining all points in the fiber \(\varphi^{-n}(\alpha)\), and \(E_{\infty}(\varphi, \alpha) := \bigcup_{n=1}^{\infty} E_{n}(\varphi, \alpha)\) is the \textit{arboreal field} associated with \(\alpha\) over \(E\).

\part{Arithmetic Inputs and Results Over Algebraic Closure}
\section{Local Inertia Orbits over Number Fields}\label{sec:inertia-orbits}

\subsection{Notations and Definitions for this Section.}

\begin{definition}
Let \(F\) be a field, let \(\phi\in F(z)\) be nonconstant, and let \(x\in\mathbb P^1(\overline F)\).
A \emph{local parameter} at \(x\in\mathbb P^{1}(\overline{F})\) is a
rational function with a simple zero at \(x\).  If \(u\) is a local
parameter at \(x\) and \(w\) one at \(\phi(x)\), then
\(w(\phi(z))=\lambda u(z)^m+\) terms of higher degree for a unique integer
\(m\geq1\) and some \(\lambda\neq0\).  We write \(\deg_x(\phi)=m\) and call
it the \emph{local degree} or \emph{multiplicity} of \(\phi\) at \(x\).
\end{definition}

Thus \(x\) is critical precisely when \(\deg_x(\phi)>1\), and
\(\mathrm{Crit}(\phi)=\{x\in\mathbb P^{1}(\overline{F}):
\deg_x(\phi)>1\}\). For any nonconstant \(\psi\), substitution of local expansions gives
\(\deg_x(\phi\circ\psi)=\deg_x(\psi)\deg_{\psi(x)}(\phi)\), and hence
\(\deg_x(\phi^r)=\prod_{i=0}^{r-1}\deg_{\phi^i(x)}(\phi)\). Since
local degree is a fiber multiplicity, we have \(1\leq\deg_x(\phi)\leq\deg(\phi)\). 

For the entirety of this section, let \(F/\mathbb Q_p\) be finite. Let \(v\) be the additive valuation normalized by \(v(F^\times)=\mathbb Z\), set \(v(0)=\infty\), and let \(|\cdot|\) be the associated absolute value. Write \(\mathcal O_F=\{x:v(x)\geq0\}\), \(\mathcal M_F=\{x:v(x)>0\}\), and \(k_F=\mathcal O_F/\mathcal M_F\); its algebraic closure is denoted by \(\overline{k_F}\). Fix an algebraic closure
\(\overline{F}\), extend \(v\) to \(\overline{F}\), and denote the
completion by \(\overline{F}_v\). Within this section,
\(\mathcal O\) and \(\mathcal M\) denote the valuation ring and maximal
ideal of \(\overline{F}_v\). For a finite Galois extension \(E/F\), write \(I(E/F)\) for its inertia group.

For \(x=[x_0:x_1]\in\mathbb P^{1}(\overline{F}_v)\), choose homogeneous
coordinates satisfying \(\max\{|x_0|,|x_1|\}=1\), and write
\(\overline{x}=[\overline{x_0}:\overline{x_1}]\in
\mathbb P^{1}(\overline{k_F})\) for its reduction.  For
\(a\in\overline{F}_v\) and \(r>0\), put
\(D(a,r)=\{z\in\overline{F}_v:|z-a|<r\}\).  The \emph{residue disk}
above \(\overline{x}\) is
\(D_{\overline{x}}=\{z\in\mathbb P^{1}(\overline{F}_v):
\overline{z}=\overline{x}\}\).  If \(\overline{x}\) is finite and
\(x_0\in\mathcal O\) is a lift, then
\(D_{\overline{x}}=\{z:|z-x_0|<1\}\); also
\(D_\infty=\{z:|z|>1\}\cup\{\infty\}\).  Translation by \(-x_0\), or
inversion when \(\overline{x}=\infty\), therefore identifies any residue
disk with \(D(0,1)\). Readers may consult
\cite[Chapter II]{SerreLocalFields} and
\cite[Section 2.1, pp. 8--10]{BenedettoAWS} for further details on local fields. We first record the following elementary local-form lemma.

\begin{lemma}[Local form under good reduction]
\label{lem:good-local-form}
Let \(\phi\in F(z)\) have good reduction, suppose
\(\overline{\phi}(0)=0\), and put \(q=\deg_0(\overline{\phi})\).  Then
\(\phi(u)=P(u)/Q(u)\) for polynomials \(P,Q\in\mathcal O_F[u]\) such that
\(Q(0)\in\mathcal O_F^\times\).  If
\(P(u)=\sum_iP_i u^i\), then \(v(P_i)>0\) for \(0\leq i<q\), while
\(v(P_q)=0\).
\end{lemma}

\begin{proof}
Choose a normalized homogeneous presentation
\(\phi=[P_h(X,Y):Q_h(X,Y)]\) with coefficients in \(\mathcal O_F\).
Good reduction means that \(\overline{P_h}\) and \(\overline{Q_h}\) have
no common projective zero and retain the degree of \(\phi\).  Since
\(0=[0:1]\) and \(\overline{\phi}(0)=0\), one has
\(\overline{P_h}(0,1)=0\) and \(\overline{Q_h}(0,1)\neq0\).  Set
\(P(u)=P_h(u,1)\) and \(Q(u)=Q_h(u,1)\).  Then \(P,Q\) are integral and
\(Q(0)\) is a unit.

Because \(\overline{Q}(0)\neq0\), the function \(\overline{Q}(u)\) is a
unit in the local ring at zero.  Hence the order of vanishing of
\(\overline{\phi}(u)=\overline{P}(u)/\overline{Q}(u)\) at zero is exactly
the order of vanishing of \(\overline{P}(u)\).  By the definition of
local degree, that order is \(q\).  Thus
\(\overline{P_i}=0\) for \(i<q\) and \(\overline{P_q}\neq0\), which is
equivalent to the asserted valuation conditions.
\end{proof}

\begin{proposition}
\label{prop:disk-bound}
Let \(F/\mathbb Q_p\) be finite, let \(\phi\in F(z)\) have good
reduction, let \(\overline{x}\in\mathbb P^{1}(\overline{k_F})\), and
let \(y\in\mathbb P^{1}(\overline{F})\) satisfy
\(\overline{y}=\overline{\phi}(\overline{x})\). Among the distinct solutions of \(\phi(z)=y\) in
\(\mathbb P^1(\overline F)\), at least one and at most
\(\deg_{\overline{x}}(\overline\phi)\) belong to
\(D_{\overline{x}}\).
\end{proposition}

\begin{proof}
Put \(m=\deg_{\overline{x}}(\overline\phi)\). After replacing \(F\)
by a finite extension, we may assume that \(y\in\mathbb P^1(F)\) and
that the finite coordinates among \(\overline{x}\) and
\(\overline{y}\) have lifts in \(\mathcal O_F\). Indeed, finite
residue-field extensions arise from finite unramified extensions, after
which we may adjoin \(y\) when it is finite. This replacement changes
neither the solutions in \(\mathbb P^1(\overline F)\) nor their
reductions. It also preserves good reduction, and scalar extension of
the reduced map does not change its local degree: a unit resultant
remains a unit, while the first nonzero exponent in a local expansion
is unchanged.

If \(\overline{x}\neq\infty\), choose
\(x_0\in\mathcal O_F\) lifting \(\overline{x}\) and put
\(\sigma(z)=z-x_0\); if \(\overline{x}=\infty\), put
\(\sigma(z)=1/z\). Define \(\tau\) similarly using a lift \(y_0\) of
\(\overline y\): take \(\tau(w)=w-y_0\) when
\(\overline y\neq\infty\), and \(\tau(w)=1/w\) otherwise. These maps
send \(\overline{x}\) and \(\overline y\), respectively, to zero, and
\(\sigma\) maps \(D_{\overline{x}}\) bijectively onto \(D(0,1)\).
Indeed, membership in the source disk becomes \(v(z-x_0)>0\) in the
finite case and \(v(1/z)>0\) in the infinite case.

Set \(\psi=\tau\circ\phi\circ\sigma^{-1}\). Integral changes of source
and target coordinates preserve good reduction and commute with
reduction; moreover, fractional-linear maps have local degree one.
Hence \(\overline\psi(0)=0,\qquad \deg_0(\overline\psi)=m\).
Lemma~\ref{lem:good-local-form} therefore gives \(\psi=P/Q\), where
\(P,Q\in\mathcal O_F[U]\), \(Q(0)\) is a unit, and, writing
\(P(U)=\sum_iP_iU^i\),
\(v(P_i)>0\quad(i<m),\qquad v(P_m)=0\).

Put \(c=\tau(y)\). Since \(y\) reduces to \(\overline y\), one has
\(v(c)>0\). If \(v(u)>0\), then every nonconstant term of \(Q(u)\)
has positive valuation, whereas \(Q(0)\) is a unit. Thus \(Q(u)\) is
a unit, and the solutions of \(\psi(u)=c\) in \(D(0,1)\) are exactly
the positive-valuation roots of \(A(U)=P(U)-cQ(U)\). All coefficients
of \(A\) are integral. For \(i<m\), both \(P_i\) and \(cQ_i\) have
positive valuation, while \(P_m\) is a unit and \(cQ_m\) has positive
valuation. Consequently, \(\operatorname{ord}_0(\overline A)=m\).

We first note that this set of roots is nonempty. Let \(H(X,Y)\) be
the homogenization of \(A(U)\). Since \(A_m\) is a unit, \(H\) is
primitive. Over a finite splitting field \(E/F\), write
\(H(X,Y)=a\prod_i(\beta_iX-\alpha_iY)\), with each pair
\((\alpha_i,\beta_i)\) primitive over \(\mathcal O_E\). Gauss's lemma
implies that \(a\) is a unit. Since
\(\overline H(0,1)=\overline A(0)=0\), some
\(\overline{\alpha_i}\) is zero. Then
\(\overline{\beta_i}\neq0\), so \([\alpha_i:\beta_i]\) reduces to
\([0:1]\). Thus \(A\) has a root in \(D(0,1)\).

We now count the positive-valuation roots of \(A\), with
multiplicity. Let \(\alpha\) be one such root and write
\(A(U)=(U-\alpha)B(U)\), where
\(A(U)=\sum_{i=0}^dA_iU^i\) and
\(B(U)=\sum_{i=0}^{d-1}B_iU^i\). Comparing coefficients from the
highest degree downward gives
\[
 B_{d-1}=A_d,\qquad B_{i-1}=A_i+\alpha B_i
 \quad(1\leq i<d).
\]
Since the \(A_i\) are integral and \(v(\alpha)>0\), induction shows
that every \(B_i\) is integral. Reducing the factorization therefore
gives \(\overline A(U)=U\overline B(U)\).

List all positive-valuation roots, with multiplicity, as
\(\alpha_1,\ldots,\alpha_s\). Applying the preceding division
successively gives
\[
 A(U)=\prod_{j=1}^s(U-\alpha_j)C(U)
 \quad\text{with }C(U)\text{ having integral coefficients},
\]
and hence \(\overline A(U)=U^s\overline C(U)\). Since
\(\operatorname{ord}_0(\overline A)=m\), it follows that \(s\leq m\).

Finally, under the bijection \(u=\sigma(z)\), the equation
\(\phi(z)=y\) for \(z\in D_{\overline{x}}\) is equivalent to
\(\psi(u)=c\) for \(u\in D(0,1)\). Its solutions therefore correspond
to the positive-valuation roots of \(A\). Their number as distinct
roots is at most \(s\leq m\), proving the proposition.
\end{proof}

\subsection{Inertia Orbits above a Periodic Reduction}
We now establish the inertia-orbit estimate used in Lemma \ref{lem:simultaneous-derangement}. Let \(F/\mathbb Q_p\) be finite, with normalized valuation \(v\). Suppose that \(\phi\in F(z)\) has good reduction, that \(c\in\mathbb P^1(F)\) is wandering, and that \(\overline c\) has exact period \(r\). Put \(q=\deg_{\overline c}(\overline\phi^{\,r})\), \(S_k(\phi,c)=\{z\in\phi^{-rk}(c):\overline z=\overline c\}\), and \(I_{rk}=I(F(\phi^{-rk}(c))/F)\). The next lemma controls the inertia orbits on the backward preimages
lying in the residue disk of \(c\). Their number stays uniformly
bounded with depth, while their cardinalities become divisible by
arbitrarily large powers of every prime dividing \(q\).

\begin{lemma}\label{lem:critical-inertia-orbits}
Let \(F/\mathbb Q_p\) be finite, let \(\phi\in F(z)\) have good reduction, and let \(c\in\mathbb P^1(F)\) be wandering. Suppose that \(\overline c\) has exact period \(r\) and that \(q\geq2\). Then the number of \(I_{rk}\)-orbits on \(S_k(\phi,c)\) is bounded independently of \(k\). For every prime \(\ell\mid q\) and every \(A\geq1\), every orbit length is divisible by \(\ell^A\) for all sufficiently large \(k\).
\end{lemma}

\begin{proof}
Choose \(\gamma\in\operatorname{PGL}_2(\mathcal O_F)\) with \(\gamma(c)=0\), and put \(h=\gamma\circ\phi^r\circ\gamma^{-1}\). Then \(h(D(0,1))\subseteq D(0,1)\), \(\overline h(0)=0\), \(h(0)\neq0\), and \(\deg_0(\overline h)=q\). By Lemma \ref{lem:good-local-form}, the denominator of an integral presentation of \(h\) is a unit on \(D(0,1)\). Its inverse has an integral power-series expansion there, and coefficient recursion gives
\[
 h(u)=\delta+\sum_{j\geq1}a_ju^j,\qquad
 0<v(\delta)<\infty,\quad v(a_j)>0\ (j<q),\quad v(a_q)=0,
\]
with \(a_j\in\mathcal O_F\) for every \(j\). Choose \(\varepsilon>0\) such that \(q\varepsilon<v(\delta)\) and \((q-j)\varepsilon<v(a_j)\) whenever \(1\leq j<q\) and \(a_j\neq0\). Then
\begin{equation}\label{eq:valuation-return}
 v(h(u))=qv(u)\qquad(0<v(u)<\varepsilon).
\end{equation}
Let \(\Delta=\{u\in\overline F_v:v(u)\geq\varepsilon\}\), a complete closed disk of radius \(R<1\). Every term of \(h(u)\) has absolute value at most \(R\) on \(\Delta\), so \(h(\Delta)\subseteq\Delta\). Moreover, factoring \(u^j-w^j\) gives \(|h(u)-h(w)|\leq\max\{|a_1|,R\}|u-w|\) for \(u,w\in\Delta\),
where \(|a_1|<1\); hence \(h\) is uniformly contracting on \(\Delta\).

There is an integer \(J\) such that no path \(0=u_0\leftarrow u_1\leftarrow\cdots\leftarrow u_J\), with \(h(u_{i+1})=u_i\), lies in \(\Delta\). Otherwise the rooted tree whose level \(n\) consists of these depth-\(n\) paths in \(\Delta\), and which has at most \(\deg(h)\) children at each vertex, would have a vertex at every level. K\"onig's infinity lemma \cite[Lemma~8.1.2, p.~245]{Diestel05} would give a path \((u_n)_{n\geq0}\) in \(\Delta\). If \(\lambda<1\) is the contraction constant, then \(|h^n(0)|=|h^n(0)-h^n(u_n)|\leq\lambda^nR\), so \(h^n(0)\to0\). Continuity applied to \(h^{n+1}(0)=h(h^n(0))\) would give \(h(0)=0\), a contradiction.

Since \(h^k\) has good reduction and
\(\overline{h^k}(0)=0\), Proposition~\ref{prop:disk-bound}, applied
with \(\overline{x}=0\) and \(y=0\), shows that
\(h^{-k}(0)\cap D(0,1)\neq\varnothing\).
Under \(\gamma\), this set corresponds to \(S_k(\phi,c)\); hence
\(S_k(\phi,c)\neq\varnothing\).
Let \(z_k\in S_k(\phi,c)\) be arbitrary, put \(u_k=\gamma(z_k)\), and set \(u_i=h^{k-i}(u_k)\). Every \(u_i\) lies in \(D(0,1)\). For \(k\geq J\), let \(j\leq J\) be the first index for which \(u_j\notin\Delta\). Forward invariance of \(\Delta\) shows that no later \(u_i\) belongs to \(\Delta\), and \eqref{eq:valuation-return} gives \(v(u_k)=v(u_j)/q^{k-j}\). The possible exit vertices lie in the finite set \(\bigcup_{i=1}^Jh^{-i}(0)\). If an exit valuation is \(A_e/B_e\) in lowest terms, then the denominator of \(v(u_k)\) is \(B_eq^{k-j}/\gcd(A_e,q^{k-j})\).
Consequently there is \(C_0\geq1\) such that this denominator is at least \(q^{k-J}/C_0\), and its \(\ell\)-adic valuation tends uniformly to infinity for every prime \(\ell\mid q\).

Let \(T_k=F(\phi^{-rk}(c))\), put \(E_k=F(z_k)\), and let \(H=\operatorname{Gal}(T_k/E_k)\). The reduced valuation denominator divides \(e(E_k/F)\). Since the inertia group of \(T_k/E_k\) is \(I_{rk}\cap H\), multiplicativity of ramification indices gives \(|I_{rk}z_k|=[I_{rk}:I_{rk}\cap H]=e(E_k/F)\). Because \(z_k\) was arbitrary, the denominator bounds apply to every inertia orbit on \(S_k(\phi,c)\).
Under \(\gamma\), the set \(S_k(\phi,c)\) is \(h^{-k}(0)\cap D(0,1)\), so Proposition \ref{prop:disk-bound} gives \(\#S_k(\phi,c)\leq q^k\). Thus the number of orbits is at most \(C_0q^J\) for \(k\geq J\); enlarging the bound handles the remaining levels. The denominator formula also proves the asserted divisibility.
\end{proof}

\section{Galois-Group Combinatorics over Number Fields}\label{sec:Group-combinatorics}

\subsection{A Frobenius Coset}
All group actions in this section are on the left. Let \(D\) be a profinite group, and let \(I\triangleleft D\) be a closed procyclic subgroup with topological generator \(\tau\). Suppose that \(\sigma\in D\) and \(Q\geq2\) satisfy \(\sigma\eta\sigma^{-1}=\eta^Q\) for every \(\eta\in I\), and let \(\mu_I\) be normalized Haar measure on \(I\). Let \(X\) be a finite set, given the discrete topology, with a continuous left \(D\)-action; let \(O\subseteq X\) be an \(I\)-orbit, and put \(s=\#O\). For \(f\geq1\), put \(\mathcal E_f(O)=\{\eta\in I:\eta\sigma^f\text{ fixes a point of }O\}\). For \(x\in O\), put \(y=\sigma^f x\) and \(E_x=\{\eta\in I:\eta y=x\}\). If \(E_x\neq\varnothing\) and \(\eta_0\in E_x\), then \(E_x=\eta_0\operatorname{Stab}_I(y)\). The stabilizer \(\operatorname{Stab}_I(y)\) is clopen because it is the inverse image of the clopen set \(\{y\}\) under the continuous orbit map \(I\to X\), \(\eta\mapsto\eta y\). Thus \(E_x\) is either empty or clopen. Since \(O\) is finite and \(\mathcal E_f(O)=\bigcup_{x\in O}E_x\), the set \(\mathcal E_f(O)\) is clopen.

The next lemma measures the fixed-point portion of the coset
\(I\sigma^f\) on a single inertia orbit. This portion is empty unless
\(\sigma^f\) preserves the orbit; when it does, its Haar measure is
\(1/\gcd(Q^f-1,\#O)\).

\begin{lemma}\label{lem:tame-coset}
Let \(D\) be a profinite group, let \(I\triangleleft D\) be closed and procyclic, and let \(\sigma\in D\) and \(Q\geq2\) satisfy \(\sigma\eta\sigma^{-1}=\eta^Q\) for every \(\eta\in I\). Let \(X\) be a finite continuous left \(D\)-set, let \(O\subseteq X\) be an \(I\)-orbit, and let \(f\geq1\). If \(\sigma^fO\neq O\), then \(\mathcal E_f(O)=\varnothing\). If \(\sigma^fO=O\), then \(\mu_I(\mathcal E_f(O))=1/\gcd(Q^f-1,s)\).
\end{lemma}

\begin{proof}
Suppose first that \(\eta\sigma^f\) fixes some \(x\in O\). Then \(\sigma^fx=\eta^{-1}x\in O\). Since \(O=I\cdot x\), we have \(\sigma^fO=\sigma^f(I\cdot x)=(\sigma^fI\sigma^{-f})\cdot\sigma^fx=I\cdot\sigma^fx\), where the last equality follows from \(I\triangleleft D\). Since \(\sigma^fx\in O\), its \(I\)-orbit is \(O\). Hence \(\sigma^fO=O\), proving that \(\mathcal E_f(O)=\varnothing\) whenever \(\sigma^fO\neq O\).

Assume now that \(\sigma^fO=O\). Choose \(x_0\in O\) and put \(H=\operatorname{Stab}_I(x_0)\). Since \(I\) is procyclic, it is abelian. Therefore \(\operatorname{Stab}_I(\xi x_0)=\xi H\xi^{-1}=H\) for every \(\xi\in I\). Thus \(H\) is the kernel of the action of \(I\) on \(O\), and \(A=I/H\) acts simply transitively on \(O\). In particular, \(A\) is cyclic of order \(s\). The image \(\overline\tau\) of \(\tau\) generates \(A\), so every point of \(O\) is uniquely of the form \(\overline\tau^{\,r}x_0\), with \(r\in\mathbb Z/s\mathbb Z\).

Since \(\sigma^fx_0\in O\), conjugation carries \(H\) onto the stabilizer of \(\sigma^fx_0\), which is again \(H\). Thus \(\sigma^fH\sigma^{-f}=\operatorname{Stab}_I(\sigma^fx_0)=H\), and conjugation by \(\sigma^f\) induces an automorphism of \(A\). Induction on \(f\) in the relation \(\sigma\eta\sigma^{-1}=\eta^Q\) gives \(\sigma^f\eta\sigma^{-f}=\eta^{Q^f}\) for every \(\eta\in I\). Hence the induced automorphism sends \(uH\) to \(u^{Q^f}H\).

There is a unique \(b\in\mathbb Z/s\mathbb Z\) such that \(\sigma^fx_0=\overline\tau^{\,b}x_0\). Consequently, \(\sigma^f(\overline\tau^{\,r}x_0)=\overline\tau^{\,Q^fr+b}x_0\). If \(\eta H=\overline\tau^{\,a}\), then \(\eta\sigma^f(\overline\tau^{\,r}x_0)=\overline\tau^{\,a+Q^fr+b}x_0\). Thus \(\eta\sigma^f\) fixes a point of \(O\) if and only if \(a+b\equiv(1-Q^f)r\pmod{s}\) for some \(r\in\mathbb Z/s\mathbb Z\).

Put \(d=\gcd(Q^f-1,s)\). Multiplication by \(1-Q^f\) on \(\mathbb Z/s\mathbb Z\) has kernel of order \(d\), and hence image of order \(s/d\). The admissible classes \(a\) form the translate of this image by \(-b\), so there are exactly \(s/d\) such classes. Thus \(\mathcal E_f(O)\) is the union of exactly \(s/d\) cosets of \(H\). The \(s\) cosets of \(H\) partition \(I\) and have equal Haar measure, so each has measure \(1/s\). Therefore \(\mu_I(\mathcal E_f(O))=(s/d)/s=1/\gcd(Q^f-1,s)\).
\end{proof}

\subsection{Derangements in Joint Preimage Fields}
For a group acting on a set, a \emph{derangement} is an element having no fixed point. The backward tree of \((\phi,c)\) has at level \(m\) the labelled set \(\{m\}\times\phi^{-m}(c)\); for \(m\geq1\), the parent of \((m,z)\) is \((m-1,\phi(z))\). Each level is finite, and every vertex has at most \(\deg(\phi)\) children. An \emph{end} of a rooted tree is an infinite ray beginning at its root; in the backward tree above, it is equivalently a sequence \((z_m)_{m\geq0}\) such that \(z_0=c\) and \(\phi(z_{m+1})=z_m\) for every \(m\). A rooted-tree automorphism fixes an end if it fixes every vertex on that ray.

K\"onig's infinity lemma states that if \(G\) is a graph and \(V_0,V_1,\ldots\) are disjoint nonempty finite subsets of \(V(G)\) such that every vertex in \(V_n\), for \(n\geq1\), has a neighbour in \(V_{n-1}\), then \(G\) contains a ray \(v_0v_1\cdots\) with \(v_n\in V_n\) for every \(n\) \cite[Lemma~8.1.2, p.~245]{Diestel05}.

\begin{lemma}\label{lem:fixed-end-criterion}
Let \(g\) be a rooted automorphism of a rooted tree with finite levels. Then \(g\) fixes an end if and only if it fixes a vertex at every level. If \(g\) fixes no end, there is an integer \(N\) such that \(g\) fixes no vertex at any level \(n\geq N\).
\end{lemma}

\begin{proof}
Only the converse requires proof. Suppose that \(g\) fixes a vertex at every level, and let \(V_n\) be the set of fixed vertices at level \(n\). The sets \(V_n\) are nonempty and finite. Since \(g\) preserves the parent map, the parent of every vertex in \(V_n\) belongs to \(V_{n-1}\). K\"onig's infinity lemma therefore gives vertices \(v_n\in V_n\) such that \(v_{n-1}\) is the parent of \(v_n\) for every \(n\geq1\). They form an end fixed by \(g\). If \(g\) fixes no end, the first assertion gives a level \(N\) without fixed vertices. No higher level has a fixed vertex, since such a vertex would have a fixed ancestor at level \(N\).
\end{proof}

Let \(F\) be a finite extension of \(\mathbb Q_p\), and put \(Q=\#k_F\). For \(1\leq j\leq u\), let \(\psi_j\in F(z)\) have degree \(1\leq d_j<p\), and let \(b_j\in\mathbb P^1(F)\). Inside a fixed algebraic closure of \(F\), put \(Z_{j,m}=\psi_j^{-m}(b_j)\), \(Z_m=\bigcup_{j=1}^uZ_{j,m}\), \(\Omega_m=F(Z_m)\), and \(\Omega=\bigcup_{m\geq0}\Omega_m\). We also need the following elementary lemma which shows that adjoining all backward preimages
simultaneously produces a Galois tower with no wild ramification.

\begin{lemma}\label{lem:joint-tower-tame}
Let \(F/\mathbb Q_p\) be finite and let \(u\geq1\). For \(1\leq j\leq u\), let \(\psi_j\in F(z)\) satisfy \(1\leq\deg(\psi_j)<p\), and let \(b_j\in\mathbb P^1(F)\). Then \(\Omega/F\) is Galois and tame.
\end{lemma}

\begin{proof}
Each \(Z_m\) is finite and Galois-stable over \(F\), so \(\Omega_m/F\) is finite Galois. Moreover, \(\psi_j(Z_{j,m})=Z_{j,m-1}\), and hence \(\Omega_{m-1}\subseteq\Omega_m\).

Put \(D_m=\operatorname{Gal}(\Omega_m/F)\). For \(m\geq1\), restriction maps \(D_m\) onto \(D_{m-1}\); let \(N_m\) be its kernel. An element of \(N_m\) fixes every point of \(Z_{j,m-1}\) and therefore permutes the points of \(\psi_j^{-1}(y)\) for each \(y\in Z_{j,m-1}\). These permutations determine the element because the sets \(\psi_j^{-1}(y)\) cover \(Z_m\), which generates \(\Omega_m\). Thus
\[
 N_m\hookrightarrow\prod_{j=1}^u\prod_{y\in Z_{j,m-1}}\operatorname{Sym}\bigl(\psi_j^{-1}(y)\bigr).
\]
Every symmetric group on the right acts on at most \(d_j<p\) points, so its order is prime to \(p\). Hence \(p\nmid|N_m|\). Since \(D_0=1\), induction in the exact sequences \(1\to N_m\to D_m\to D_{m-1}\to1\) gives \(p\nmid|D_m|\) for every \(m\).

Absolute wild inertia is pro-\(p\) by \cite[Proposition~7.5.1, p.~409]{NSW08}. Its image in each \(D_m\) is therefore trivial. Thus every \(\Omega_m/F\) is tame, and so is their union \(\Omega/F\). The union is Galois because the extensions \(\Omega_m/F\) are nested and Galois.
\end{proof}

Put \(D=\operatorname{Gal}(\Omega/F)\), let \(I\) be the subgroup of \(D\) acting trivially on the residue field of \(\Omega\), and let \(\mu_I\) be normalized Haar measure on \(I\). The next lemma packages the tame local Galois structure of the tower: its inertia group is procyclic, surjects onto inertia at every finite
level, and admits a Frobenius lift whose conjugation action is the
\(Q\)-power map.

\begin{lemma}\label{lem:tame-tower-inertia}
Let \(F/\mathbb Q_p\) be finite with \(\#k_F=Q\), and let \(u\geq1\). For \(1\leq j\leq u\), let \(\psi_j\in F(z)\) satisfy \(1\leq\deg(\psi_j)<p\), and let \(b_j\in\mathbb P^1(F)\). Restriction maps \(I\) onto \(I(T/F)\) for every finite Galois field \(T\) with \(F\subseteq T\subseteq\Omega\). Moreover, \(I\) is procyclic, and there exist a topological generator \(\tau\) of \(I\) and a Frobenius lift \(\sigma\in D\) such that \(\sigma\eta\sigma^{-1}=\eta^Q\) for every \(\eta\in I\).
\end{lemma}

\begin{proof}
Let \(F^{\mathrm t}\) be the maximal tamely ramified extension of \(F\) in a fixed algebraic closure, let \(F^{\mathrm{nr}}\subseteq F^{\mathrm t}\) be its maximal unramified subextension, and put \(I_F^{\mathrm t}=\operatorname{Gal}(F^{\mathrm t}/F^{\mathrm{nr}})\). Since \(\Omega/F\) is tame, restriction gives \(\rho:\operatorname{Gal}(F^{\mathrm t}/F)\to D\).

Let \(T/F\) be a finite Galois intermediate extension. The image of \(I_F^{\mathrm t}\) in \(\operatorname{Gal}(T/F)\) fixes exactly \(T\cap F^{\mathrm{nr}}\), the maximal unramified subextension of \(T/F\), and is therefore \(I(T/F)\). In particular, its order is the ramification index \(e(T/F)\), as also follows from \cite[Chapter~I, Section~7, Corollary to Proposition~21, p.~22]{SerreLocalFields}.

Put \(J=\rho(I_F^{\mathrm t})\). Since tame inertia acts trivially on residues, \(J\subseteq I\). For \(g\in I\) and every finite Galois intermediate field \(T\), the preceding paragraph gives an element of \(J\) having the same restriction to \(T\) as \(g\). The kernels of restriction to such fields form a neighbourhood basis in \(D\), so \(J\) is dense in \(I\). It is also compact and therefore closed. Hence \(J=I\), and restriction maps \(I\) onto \(I(T/F)\) for every \(T\).

By \cite[Proposition~7.5.2 and Theorem~7.5.3, pp.~410--411]{NSW08}, \(I_F^{\mathrm t}\) has a topological generator \(\tau_F\), and \(\operatorname{Gal}(F^{\mathrm t}/F)\) contains a lift \(\widetilde\sigma\) of the residue automorphism \(x\mapsto x^Q\) satisfying \(\widetilde\sigma\tau_F\widetilde\sigma^{-1}=\tau_F^Q\). Put \(\tau=\rho(\tau_F)\) and \(\sigma=\rho(\widetilde\sigma)\). Then \(\tau\) topologically generates \(I\), and \(\sigma\) induces \(x\mapsto x^Q\) on the residue field of \(\Omega\). The required identity holds on the dense cyclic subgroup generated by \(\tau\), and hence on all of \(I\) by continuity.
\end{proof}

Fix a Frobenius lift \(\sigma\in D\) such that \(\sigma\eta\sigma^{-1}=\eta^Q\) for every \(\eta\in I\). Let \(\phi\in F(z)\) have good reduction, let \(c\in\mathbb P^1(F)\), and suppose that \(F(\phi^{-n}(c))\subseteq\Omega\) for every \(n\geq0\). For \(f\geq1\), put \(Y_n=\phi^{-n}(c)\), \(\mathcal B_n=\{\eta\in I:\eta\sigma^f\text{ fixes a point of }Y_n\}\), and \(\mathcal B=\bigcap_{n\geq0}\mathcal B_n\). If \(\overline c\) is periodic, let \(r\) be its exact period, put \(S_k=\{z\in Y_{rk}:\overline z=\overline c\}\) and \(T_k=F(Y_{rk})\), and suppose that there are a prime \(\ell\) and integers \(A,C\geq1\) such that, for every sufficiently large \(k\), the group \(I(T_k/F)\) has at most \(C\) orbits on \(S_k\), each of length divisible by \(\ell^A\). Assume also that \(Q^f\equiv1\pmod{\ell^A}\). The next lemma measures the elements of the Frobenius coset that fix
an infinite backward branch. There are none when \(\overline c\) is
nonperiodic; in the periodic case they form a closed set of Haar
measure at most \(C/\ell^A\).

\begin{lemma}\label{lem:fixed-end-measure}
Let \(F,\Omega,I,\sigma,Q\) be as in Lemma \ref{lem:tame-tower-inertia}. Let \(\phi\in F(z)\) have good reduction, let \(c\in\mathbb P^1(F)\), suppose that \(F(\phi^{-n}(c))\subseteq\Omega\) for every \(n\geq0\), and let \(f\geq1\). Then \(\mathcal B\) is closed. If \(\overline c\) is not periodic, then \(\mathcal B=\varnothing\). Suppose that \(\overline c\) has exact period \(r\), that \(I(T_k/F)\) has at most \(C\) orbits on \(S_k\), each of length divisible by \(\ell^A\), for every sufficiently large \(k\), and that \(Q^f\equiv1\pmod{\ell^A}\). Then \(\mu_I(\mathcal B)\leq C/\ell^A\).
\end{lemma}

\begin{proof}
Every element of \(D\) fixes \(c\) and commutes with \(\phi\), so it acts as a rooted automorphism of the backward tree of \(c\). For every \(n\), the fixed-point condition defining \(\mathcal B_n\) depends only on the action of \(I\) on the finite set \(Y_n\), so \(\mathcal B_n\) is clopen. Hence \(\mathcal B\) is closed. An element \(\eta\) belongs to \(\mathcal B\) precisely when \(\eta\sigma^f\) fixes a vertex at every level. By Lemma \ref{lem:fixed-end-criterion}, this is equivalent to fixing an end.

Suppose that \(\eta\in\mathcal B\), put \(\gamma=\eta\sigma^f\), and choose a fixed end \((z_n)_{n\geq0}\). Since \(\eta\) acts trivially on residues and \(\sigma^f\) induces the \(Q^f\)-power Frobenius, the equality \(\gamma z_n=z_n\) implies \(\overline z_n\in\mathbb P^1(\mathbb F_{Q^f})\). Good reduction gives \(\overline\phi(\overline z_{n+1})=\overline z_n\).

Fix \(n\). The sequence \(\overline z_n,\overline z_{n+1},\ldots\) lies in the finite set \(\mathbb P^1(\mathbb F_{Q^f})\), so there are \(n\leq a<b\) with \(\overline z_a=\overline z_b\). Hence \(\overline\phi^{\,b-a}(\overline z_b)=\overline z_b\), and \(\overline z_n=\overline\phi^{\,b-n}(\overline z_b)\) lies on the same cycle. Since \(\overline\phi^{\,n}(\overline z_n)=\overline c\), this is the cycle containing \(\overline c\). Thus \(\overline c\) is periodic, proving that \(\mathcal B=\varnothing\) when \(\overline c\) is not periodic. In the periodic case, \(\overline\phi^{\,rk}\) is the identity on this cycle. Since \(\overline\phi^{\,rk}(\overline z_{rk})=\overline c\), we obtain \(\overline z_{rk}=\overline c\), and hence \(z_{rk}\in S_k\), for every \(k\geq1\).

Choose \(k\) for which the orbit hypotheses hold. The set \(Y_{rk}\) is Galois-stable, so \(T_k/F\) is finite Galois; by assumption, \(T_k\subseteq\Omega\). The set \(S_k\) is \(I\)-stable because inertia acts trivially on residues. By Lemma \ref{lem:tame-tower-inertia}, restriction maps \(I\) onto \(I(T_k/F)\), so the \(I\)-orbits on \(S_k\) are exactly the \(I(T_k/F)\)-orbits. Thus there are at most \(C\) such orbits, and every one has cardinality divisible by \(\ell^A\).

If \(\eta\in\mathcal B\), the level-\(rk\) vertex of an \(\eta\sigma^f\)-fixed end lies in \(S_k\) and is fixed by \(\eta\sigma^f\). Therefore
\[
 \mathcal B\subseteq\bigcup_{O\in I\backslash S_k}\mathcal E_f(O),
\]
where each \(O\) is regarded as an \(I\)-orbit in the finite \(D\)-set \(Y_{rk}\). If \(\sigma^fO\neq O\), Lemma \ref{lem:tame-coset} gives \(\mathcal E_f(O)=\varnothing\). If \(\sigma^fO=O\), it gives \(\mu_I(\mathcal E_f(O))=1/\gcd(Q^f-1,\#O)\). Both integers in the greatest common divisor are divisible by \(\ell^A\), so \(\mu_I(\mathcal E_f(O))\leq\ell^{-A}\). The union contains at most \(C\) orbits, and hence \(\mu_I(\mathcal B)\leq C/\ell^A\).
\end{proof}

Let \(L\) be a number field and let \(s\geq1\). For \(1\leq j\leq s\), let \(\phi_j\in L(z)\) have degree at least \(2\), and let \(c_j\in\mathbb P^1(L)\) be a wandering critical point of \(\phi_j\). Put \(\mathcal L_0=L\) and \(Y_{j,0}=\{c_j\}\). For \(m\geq1\), let \(\mathcal L_m\) be the compositum of the fields \(L_m(\phi_j,c_j)\), put \(Y_{j,m}=\phi_j^{-m}(c_j)\), and put \(\mathcal X_m=\bigcup_{j=1}^sY_{j,m}\). The following final lemma assembles the preceding local estimates into a single
global derangement. At every sufficiently deep level, one Galois
element simultaneously moves every preimage in all the chosen backward
orbits.

\begin{lemma}\label{lem:simultaneous-derangement}
Let \(L\) be a number field and let \(s\geq1\). For \(1\leq j\leq s\), let \(\phi_j\in L(z)\) have degree at least \(2\), and let \(c_j\in\mathbb P^1(L)\) be a wandering critical point of \(\phi_j\). Then, for every sufficiently large \(m\), \(\operatorname{Gal}(\mathcal L_m/L)\) contains a derangement of \(\mathcal X_m\).
\end{lemma}

\begin{proof}
Put \(d_j=\deg(\phi_j)\) and \(e_j=\deg_{c_j}(\phi_j)\), so \(e_j\geq2\). Choose homogeneous presentations of the \(\phi_j\) and linear-fractional local coordinates at \(c_j\) and \(\phi_j(c_j)\), all over \(L\). Exclude the finitely many finite places at which one of these presentations or coordinate changes is not integral with unit resultant or determinant, or the coefficient of the first nonzero local term of some \(\phi_j\) at \(c_j\) is not a unit. At every remaining place, each \(\phi_j\) has good reduction and \(\deg_{\overline c_j}(\overline\phi_j)=e_j\). Exclude also the places of residue characteristic at most \(\max_jd_j\), and choose a finite place \(v\) outside the resulting finite set. Put \(F=L_v\), let \(p\) be its residue characteristic, and put \(Q=\#k_v\). Fix an embedding \(\iota:\overline L\hookrightarrow\overline F\) extending \(L\hookrightarrow F\), and use bars for reduction with respect to the induced valuation.

Let \(P\) be the set of indices \(j\) for which \(\overline c_j\) is periodic under \(\overline\phi_j\). For \(j\in P\), let \(r_j\) be its exact period and put \(q_j=\deg_{\overline c_j}(\overline\phi_j^{\,r_j})\). By multiplicativity of local degrees, \(q_j=\prod_{a=0}^{r_j-1}\deg_{\overline\phi_j^{\,a}(\overline c_j)}(\overline\phi_j)\). The factor with \(a=0\) is \(e_j\), so \(q_j\geq2\). Every factor is at most \(d_j<p\), and therefore \(p\nmid q_j\). For \(k\geq1\), put \(S_{j,k}=\{z\in\phi_j^{-r_jk}(c_j):\overline z=\overline c_j\}\). Lemma \ref{lem:critical-inertia-orbits} gives an integer \(C_j\) bounding, independently of \(k\), the number of \(I(F(\phi_j^{-r_jk}(c_j))/F)\)-orbits on \(S_{j,k}\). Choose a prime \(\ell_j\mid q_j\) and an integer \(A_j\geq1\) for each \(j\in P\) so that \(\sum_{j\in P}C_j/\ell_j^{A_j}<1\). For these choices, the same lemma shows that every such inertia orbit has length divisible by \(\ell_j^{A_j}\) for every sufficiently large \(k\). Since \(\ell_j\neq p\), the number \(Q\) is a unit modulo \(\ell_j^{A_j}\). Choose \(f\geq1\) divisible by the multiplicative order of \(Q\) modulo \(\ell_j^{A_j}\) for every \(j\in P\); if \(P=\varnothing\), take \(f=1\).

The set \(\mathcal X_m\) is Galois-stable, so \(\mathcal L_m=L(\mathcal X_m)\) is finite Galois over \(L\). Moreover, \(\phi_j(Y_{j,m})=Y_{j,m-1}\), and hence \(\mathcal L_{m-1}\subseteq\mathcal L_m\). Put \(\Omega_m=F\iota(\mathcal L_m)=F(\iota(\mathcal X_m))\) and \(\Omega=\bigcup_{m\geq0}\Omega_m\). The embedding \(\iota\) determines compatible places \(w_m\) of \(\mathcal L_m\) above \(v\). The closure of \(\iota(\mathcal L_m)\) in \(\overline F\) is \(F\iota(\mathcal L_m)\): it contains the closure \(F\) of \(L\), while \(F\iota(\mathcal L_m)\) is finite over \(F\) and hence closed. Thus \(\Omega_m\) is canonically the completion \((\mathcal L_m)_{w_m}\). By \cite[Chapter~II, Section~9, Proposition~(9.6), p.~170]{NeukirchANT}, restriction identifies \(\operatorname{Gal}(\Omega_m/F)\) with the decomposition subgroup \(D(w_m/v)\subseteq\operatorname{Gal}(\mathcal L_m/L)\), compatibly with restriction and with their actions on \(\iota(\mathcal X_m)\). We henceforth identify \(Y_{j,m}\) and \(\mathcal X_m\) with their images under \(\iota\).

Lemma \ref{lem:joint-tower-tame} shows that \(\Omega/F\) is tame. Put \(D=\operatorname{Gal}(\Omega/F)\), let \(I\) be its inertia subgroup, and let \(\mu_I\) be normalized Haar measure on \(I\). Lemma \ref{lem:tame-tower-inertia} gives a topological generator \(\tau\) of \(I\), a Frobenius lift \(\sigma\in D\) satisfying \(\sigma\eta\sigma^{-1}=\eta^Q\) for every \(\eta\in I\), and surjective restriction from \(I\) to the inertia group of every finite Galois intermediate extension of \(\Omega/F\).

For \(1\leq j\leq s\) and \(n\geq0\), put \(\mathcal B_{j,n}=\{\eta\in I:\eta\sigma^f\text{ fixes a point of }Y_{j,n}\}\) and \(\mathcal B_j=\bigcap_{n\geq0}\mathcal B_{j,n}\). Lemmas \ref{lem:critical-inertia-orbits} and \ref{lem:fixed-end-measure} show that \(\mathcal B_j\) is closed, that \(\mathcal B_j=\varnothing\) for \(j\notin P\), and that \(\mu_I(\mathcal B_j)\leq C_j/\ell_j^{A_j}\) for \(j\in P\). Consequently \(\mu_I(\bigcup_{j=1}^s\mathcal B_j)\leq\sum_{j\in P}C_j/\ell_j^{A_j}<1\).

Choose \(\eta\in I\setminus\bigcup_j\mathcal B_j\) and put \(\gamma=\eta\sigma^f\). For every \(j\), the condition \(\eta\notin\mathcal B_j\) and Lemma \ref{lem:fixed-end-criterion} show that \(\gamma\) fixes no end of the backward tree of \(c_j\). The same lemma therefore gives an integer \(m_j\) such that \(\gamma\) fixes no vertex of that tree at any level \(m\geq m_j\). Put \(m_0=\max_jm_j\). For \(m\geq m_0\), the restriction \(\gamma_m=\gamma|_{\Omega_m}\) fixes no point of \(\mathcal X_m\). Under the decomposition-group identification above, \(\gamma_m\) corresponds to an element of \(D(w_m/v)\subseteq\operatorname{Gal}(\mathcal L_m/L)\) with the same action on \(\mathcal X_m\). This element is therefore a derangement of \(\mathcal X_m\).
\end{proof}

\section{Proof of the Simultaneous Hasse Principle over Number Fields}\label{sec:proof-Hasse-principle}
Let \(K_0\) denote the field in the statement of Theorem \ref{thm:simultaneousHasse}. Choose a finite Galois extension \(K/K_0\) containing every critical point in the sets \(\mathcal C_i\). It is enough to construct the required primes over \(K\). Indeed, from any positive-density set of primes of \(K\), discard the primes of relative residue degree greater than one over \(K_0\) and contract the rest to \(K_0\). The discarded primes have density zero, and contraction has bounded fibers, so the contracted set still has positive density. If a critical point \(c\) were periodic at a prime \(\mathfrak q\) of \(K_0(c)\) above one of the contracted primes, choose a prime \(\mathfrak q_K\) of \(K\) above \(\mathfrak q\). Galois transitivity supplies \(\sigma\in\operatorname{Gal}(K/K_0)\) carrying the retained prime above \(K_0\) to \(\mathfrak q_K\). Then \(\sigma^{-1}(c)\), which is again a nonperiodic critical point of the same map, would be periodic at the retained prime, a contradiction. We may therefore work over \(K\), where all the points in every \(\mathcal C_i\) are rational.

First consider a strictly preperiodic point \(c\in\mathcal C_i\). Let \(a\geq1\) be its tail length and \(b\geq1\) the period of its eventual cycle. Then \(\varphi_i^{a+b}(c)=\varphi_i^a(c)\), whereas \(\varphi_i^b(c)\neq c\). At a prime of good reduction, if \(\overline c\) were periodic, then \(\overline\varphi_i^a\) would be injective on its cycle; reducing the equality would give \(\overline\varphi_i^b(\overline c)=\overline c\). Since \(\varphi_i^b(c)\) and \(c\) are distinct, this can occur at only finitely many primes. Thus every strictly preperiodic point in the sets \(\mathcal C_i\) is nonperiodic modulo all but finitely many primes.

Let \(\Gamma_{\mathrm w}=\{(i,c):c\in\mathcal C_i\text{ is wandering}\}\). If \(\Gamma_{\mathrm w}=\varnothing\), the preceding paragraph proves the theorem. Otherwise, for \(m\geq1\), let \(L_m\) be the compositum of the fields \(K_m(\varphi_i,c)\), for \((i,c)\in\Gamma_{\mathrm w}\), and put \(X_m=\bigcup_{(i,c)\in\Gamma_{\mathrm w}}\varphi_i^{-m}(c)\). Choose \(m\) sufficiently large for Lemma \ref{lem:simultaneous-derangement}. Applied over \(K\) to all pairs indexed by \(\Gamma_{\mathrm w}\), the lemma gives a derangement \(\delta\in\operatorname{Gal}(L_m/K)\) of \(X_m\). Since \(X_m\) is Galois-stable, every conjugate of \(\delta\) is also a derangement: if \(a\delta a^{-1}\) fixed \(x\in X_m\), then \(\delta\) would fix \(a^{-1}x\).

Let \(\mathcal D\) be the conjugacy class of \(\delta\). Chebotarev's theorem gives a positive-density set \(\mathcal S_m\) of primes \(\mathfrak p\) of \(K\) with Frobenius class \(\mathcal D\); see \cite[Theorem~6.3.1]{FrJar} or \cite{SL96}. Remove the finitely many primes of bad reduction and those at which two distinct points of \(X_m\) have the same reduction. For every remaining \(\mathfrak p\in\mathcal S_m\), the Frobenius--factorization criterion shows that none of the equations \(\varphi_i^m(x)=c\), with \((i,c)\in\Gamma_{\mathrm w}\), has a solution in \(\mathbb P^1(k_{\mathfrak p})\). If \(\overline c\) were periodic, its cycle would contain an \(m\)-th predecessor of \(\overline c\), contradicting this conclusion. Hence every wandering point in the sets \(\mathcal C_i\) is nonperiodic modulo every such \(\mathfrak p\).

Finally, remove the finitely many exceptional primes for the strictly preperiodic points. The remaining set still has positive density and consists of primes of common good reduction at which no point of any \(\mathcal C_i\) is periodic. \(\square\)

\section{Critical-Residue Avoidance over Number Fields} \label{sec:Hasse-to-DML}
We first prove the uniform critical-residue avoidance needed for Theorem \ref{thm:padicattracting}.
For a geometric point \(\alpha\), the congruence below is understood at every prime of \(K(\alpha)\) above \(\mathfrak p\).

\begin{theorem}\label{thm:intermediate}
Let \(K\) be a number field, let \(\varphi_{1},\ldots,\varphi_{g}:\mathbb{P}_{K}^{1}\longrightarrow\mathbb{P}_{K}^ {1}\) be rational maps of degree at least 2, and for every \(1 \leq i \leq g\), let \(\mathcal{C}_{i}\) be the set of all the critical points of \(\varphi_{i}\). Let \({\mathcal{T}}_{1},\ldots,{\mathcal{T}}_{g}\) be finite subsets of \(\mathbb{P}^{1}(K)\) such that no \({\mathcal{T}}_{i}\) contains any \(\varphi_{i}\)-preperiodic points. Then, there is a positive integer \(M\) and a set of primes \({\mathcal{P}}\) of \(K\) with positive density such that for any \(i=1,\ldots,g\), any \(\gamma\in{\mathcal{T}}_{i}\), any \(\alpha\in{\mathcal{C}}_{i}\), any \(\mathfrak{p}\in{\mathcal{P}}\), and any \(m\geq M\),
\[
\varphi_{i}^{m}(\gamma)\not\equiv\alpha\pmod{\mathfrak{p}}.
\]
\end{theorem}

\subsection{Proof of Theorem \ref{thm:intermediate}}

\subsubsection{Reduction to Finite Extensions.}
Let $K_0$ denote the original ground field. Choose a finite Galois
extension $K/K_0$ containing every critical point. We may work over
$K$ and over any further finite extension $L/K$. Indeed, from a
positive-density set of primes of $L$, discard the primes of relative
residue degree greater than one over $K_0$ and contract the remainder
to $K_0$. The discarded primes have density zero, and contraction has
bounded fibers, so the resulting set still has positive density. Let
$\mathfrak r$ be a retained prime, put
$\mathfrak p=\mathfrak r\cap\mathfrak{o}_{K_0}$ and
$\mathfrak p_K=\mathfrak r\cap\mathfrak{o}_K$, and suppose that
$\varphi_i^m(\gamma)\equiv\alpha$ at a prime $\mathfrak q$ of
$K_0(\alpha)$ above $\mathfrak p$. Choose a prime $\mathfrak q_K$ of
$K$ above $\mathfrak q$ and
$\sigma\in\operatorname{Gal}(K/K_0)$ with
$\sigma(\mathfrak p_K)=\mathfrak q_K$. Applying $\sigma^{-1}$ gives
$\varphi_i^m(\gamma)\equiv\sigma^{-1}(\alpha)$ at $\mathfrak p_K$,
hence at $\mathfrak r$, contradicting avoidance over $L$. Thus it
suffices to work over $K$, and henceforth every critical point is
$K$-rational.

\subsubsection{Reduction to Fixed Periodic Critical Points.}
For each $i$, let $b_i$ be the least common multiple of the exact
periods of the $\varphi_i$-periodic critical points, with $b_i=1$ if
there are none, and put $\psi_i=\varphi_i^{b_i}$,
$\mathcal{A}_i=\operatorname{Crit}(\varphi_i)$, and
$\mathcal{T}_i'=\bigcup_{0\leq r<b_i}\varphi_i^r(\mathcal{T}_i)$.
Every point of $\mathcal{T}_i'$ is $\psi_i$-wandering, every periodic
point of $\mathcal{A}_i$ is fixed by $\psi_i$, and every wandering
point of $\mathcal{A}_i$ is a critical point of $\psi_i$.

It suffices to prove the corresponding assertion for
$(\psi_i,\mathcal{A}_i,\mathcal{T}_i')$. Indeed, let
$B=\max_i b_i$, suppose that the assertion for these data has bound
$N$, and write $m=b_i n+r$ with $0\leq r<b_i$. If $m\geq BN$, then
$n\geq N$ and
$\varphi_i^m(\gamma)=\psi_i^n(\varphi_i^r(\gamma))$, where
$\varphi_i^r(\gamma)\in\mathcal{T}_i'$. We therefore replace
$(\varphi_i,\mathcal{C}_i,\mathcal{T}_i)$ by
$(\psi_i,\mathcal{A}_i,\mathcal{T}_i')$ and retain the notation
$(\varphi_i,\mathcal{C}_i,\mathcal{T}_i)$ for the replacements.

Thus every point of $\mathcal{T}_i$ is $\varphi_i$-wandering, every
periodic point of $\mathcal{C}_i$ is fixed by $\varphi_i$, and every
wandering point of $\mathcal{C}_i$ is a critical point of
$\varphi_i$. By Theorem \ref{thm:simultaneousHasse}, the wandering points of all the
$\mathcal{C}_i$ remain nonperiodic modulo every prime in a
positive-density set. Every strictly preperiodic point of
$\mathcal{C}_i$ remains nonperiodic outside a finite set of primes.
Delete also the primes of bad reduction, those of residue
characteristic at most $\max_i\deg(\varphi_i)$, and those at which the
local degree at a fixed critical point changes.

\subsubsection{The Ramification Proposition.}
We use the following proposition of Benedetto, Ghioca, Hutz, Kurlberg,
Scanlon, and Tucker \cite[Proposition~3.4]{BGHKST13}.

\noindent\textbf{Proposition 2.}
Let $K$ be a number field, let $\widetilde{\mathfrak p}$ be a prime of
$\mathfrak{o}_{\overline K}$, and let
$\varphi\colon\mathbb P^1\to\mathbb P^1$ be a rational map defined
over $K$, with good reduction at
$\mathfrak p=\widetilde{\mathfrak p}\cap\mathfrak{o}_K$, such that
$2\leq\deg(\varphi)<\operatorname{char}(k_{\mathfrak p})$. Let
$\mathcal{A}\subseteq\mathbb P^1(K)$ be finite and suppose:
\begin{itemize}
\item every nonperiodic point of $\mathcal{A}$ remains nonperiodic
modulo $\mathfrak p$;
\item every periodic point of $\mathcal{A}$ is fixed by $\varphi$, and
its ramification index for $\varphi$ is unchanged by reduction.
\end{itemize}
Then there is a finite extension $E/K$ such that, for every finite
extension $L/E$, there is an integer $M\geq1$ with the following
property. Whenever $m\geq M$ and $\beta\in\mathbb P^1(\overline K)$
satisfy $\varphi^m(\beta)\in\mathcal{A}$ and
$\varphi^t(\beta)\notin\mathcal{A}$ for $0\leq t<m$, the prime
$\mathfrak r=\widetilde{\mathfrak p}\cap\mathfrak{o}_{L(\beta)}$ is
unramified over
$\mathfrak q=\widetilde{\mathfrak p}\cap\mathfrak{o}_L$, and
$[\mathfrak{o}_{L(\beta)}/\mathfrak r:
\mathfrak{o}_L/\mathfrak q]>1$.

Fix a remaining prime $\mathfrak p$ of $K$ and a prime
$\widetilde{\mathfrak p}$ of $\mathfrak{o}_{\overline K}$ above it.
The proposition applies to every pair $(\varphi_i,\mathcal{C}_i)$.
Let $E_i/K$ be the resulting extensions, let $L$ be their compositum,
and put $\mathfrak q=\widetilde{\mathfrak p}\cap\mathfrak{o}_L$. For
each $i$, choose the corresponding integer $M_i$ for $L/E_i$, and set
$M=\max_i M_i$. Hence, whenever
$\varphi_i^M(\beta)\in\mathcal{C}_i$ and
$\varphi_i^t(\beta)\notin\mathcal{C}_i$ for $0\leq t<M$, the prime
$\mathfrak s_\beta=
\widetilde{\mathfrak p}\cap\mathfrak{o}_{L(\beta)}$ is unramified
over $\mathfrak q$ and
$[k_{\mathfrak s_\beta}:k_{\mathfrak q}]>1$.

We now collect all depth-\(M\) first-entry points into one finite
Galois-stable set \(Y\).
We choose a Frobenius class acting without fixed points on \(Y\), and
then show that any sufficiently late critical residue hit would produce
a Frobenius-fixed point of \(Y\).

\subsubsection{The First-Entry Field and Chebotarev.}
Let $Y$ be the union, over all $i$, of the points
$\beta\in\mathbb P^1(\overline L)$ satisfying
$\varphi_i^M(\beta)\in\mathcal{C}_i$ and
$\varphi_i^t(\beta)\notin\mathcal{C}_i$ for $0\leq t<M$, and put
$F=L(Y)$. The set $Y$ is finite because each $\varphi_i^M$ has finite
fibers. For $\sigma\in\operatorname{Gal}(\overline L/L)$, one has
$\varphi_i^t(\sigma\beta)=\sigma(\varphi_i^t(\beta))$, and $\sigma$
fixes $\mathcal{C}_i$ pointwise; hence $\sigma(Y)=Y$. Thus $F/L$ is
Galois.

Put $\mathfrak Q=\widetilde{\mathfrak p}\cap\mathfrak{o}_F$. The
completion $F_{\mathfrak Q}$ is the compositum of the extensions
$L(\beta)_{\mathfrak s_\beta}/L_{\mathfrak q}$, each of which is
unramified. Consequently $F/L$ is unramified at $\mathfrak q$. Let
$\tau=\operatorname{Frob}_{\mathfrak Q}\in\operatorname{Gal}(F/L)$.
If $\tau$ fixed some $\beta\in Y$, then its restriction to
$L(\beta)$ would be the identity, whereas its residue action is the
Frobenius of $k_{\mathfrak s_\beta}/k_{\mathfrak q}$, whose order is
$[k_{\mathfrak s_\beta}:k_{\mathfrak q}]>1$. This is impossible.
Since $Y$ is Galois-stable, no conjugate of $\tau$ fixes a point of
$Y$.

Chebotarev's theorem gives a positive-density set $\mathcal{S}$ of
primes $\mathfrak r$ of $L$ having Frobenius conjugacy class $[\tau]$.
Delete the primes of bad reduction for any $\varphi_i$ and the
finitely many primes at which reduction is not injective on $Y$.

\subsubsection{The First-Entry Claim.}
\noindent\textbf{Claim 1.}
Suppose that $z\in\mathbb P^1(L)$ and $\varphi_i^n(z)$ is congruent
modulo some $\mathfrak r\in\mathcal{S}$ to a point of
$\mathcal{C}_i$ for some $n\geq0$. Then $\varphi_i^t(z)$ is congruent
modulo $\mathfrak r$ to a point of $\mathcal{C}_i$ for some
$0\leq t<M$.

\noindent\textit{Proof of Claim 1.}
Choose the least such $n$. There is nothing to prove if $n<M$.
Otherwise put $w=\varphi_i^{n-M}(z)$. Minimality gives
$\varphi_i^t(w)\not\equiv c\pmod{\mathfrak r}$ for every
$c\in\mathcal{C}_i$ and $0\leq t<M$. Choose a prime $\mathfrak R$ of
$F$ above $\mathfrak r$ and let
$\tau_{\mathfrak r}=\operatorname{Frob}_{\mathfrak R}$. Since
$\mathfrak r\in\mathcal{S}$, the element $\tau_{\mathfrak r}$ belongs
to the conjugacy class $[\tau]$. Put $x_j=\overline{\varphi_i^j(w)}$ for $0\leq j\leq M$, and choose
$c\in\mathcal{C}_i$ with $x_M=\overline c$. Starting with
$\beta_M=c$, Proposition  \ref{prop:disk-bound} applied recursively for
$j=M-1,\ldots,0$ gives $\beta_j$ reducing to $x_j$ and satisfying
$\varphi_i(\beta_j)=\beta_{j+1}$. Thus $\beta=\beta_0$ reduces to
$w$, satisfies $\varphi_i^M(\beta)=c$, and has
$\varphi_i^j(\beta)$ reducing to $\varphi_i^j(w)$ for $0\leq j<M$.
The absence of earlier congruences with $\mathcal{C}_i$ therefore
gives $\varphi_i^j(\beta)\notin\mathcal{C}_i$ for $0\leq j<M$.
Hence $\beta\in Y$. The reduction of
$\beta$ is therefore $\overline w\in\mathbb P^1(k_{\mathfrak r})$.
On the residue field at $\mathfrak R$, the automorphism
$\tau_{\mathfrak r}$ acts as the $\lvert k_{\mathfrak r} \rvert$-power map, which
fixes $\mathbb P^1(k_{\mathfrak r})$ pointwise. Hence
$\tau_{\mathfrak r}(\beta)$ also reduces to $\overline w$. Moreover,
$\tau_{\mathfrak r}(\beta)\in Y$ because $Y$ is Galois-stable.
Injectivity of reduction on $Y$ now gives
$\tau_{\mathfrak r}(\beta)=\beta$, contradicting the fact that no
element of $[\tau]$ fixes a point of $Y$. This proves the claim.

\subsubsection{Exclusion of Long Critical Hits.}
Let $\mathcal{W}_i\subseteq\mathcal{C}_i$ be the set of
$\varphi_i$-wandering points and put $h=\max \{\lvert \mathcal{W}_i \rvert: 1 \leq i \leq g\}$.
Remove from $\mathcal{S}$ the finitely many primes satisfying either
of the following:
\begin{enumerate}
\item[(a)] $\varphi_i^t(\delta)\equiv c$ for some
$\delta\in\mathcal{T}_i$, some preperiodic $c\in\mathcal{C}_i$, and
some $0\leq t<M$;
\item[(b)] $\varphi_i^s(a)\equiv b$ for some $a\in\mathcal{W}_i$,
$b\in\mathcal{C}_i$, and $1\leq s\leq M$ such that
$\varphi_i^s(a)\neq b$ over $L$.
\end{enumerate}

Only finitely many primes are removed. In (a),
$\varphi_i^t(\delta)$ is wandering, whereas $c$ is preperiodic, so
the two points are distinct. In (b), $\varphi_i^s(a)$ and $b$ are
distinct by hypothesis. Two distinct points of $\mathbb P^1(L)$ have
the same reduction at only finitely many primes of $L$; since all
indexing sets and exponent ranges are finite, the union of exceptional
primes is finite.

Suppose, at a remaining prime, that
$\varphi_i^n(\delta)\equiv\alpha$ for some
$\delta\in\mathcal{T}_i$, $\alpha\in\mathcal{C}_i$, and
$n\geq(h+1)M$. Claim 1 gives $0\leq t_0<M$ and
$a_0\in\mathcal{C}_i$ such that
$\varphi_i^{t_0}(\delta)\equiv a_0$. By (a), $a_0$ is wandering; if
$\mathcal{W}_i=\varnothing$, this is already a contradiction. Put
$R_0=n-t_0$. Then $R_0\geq hM+1$ and
$\varphi_i^{R_0}(a_0)\equiv\alpha$.

Suppose that $a_j\in\mathcal{W}_i$,
$R_j\geq(h-j)M+1$, and
$\varphi_i^{R_j}(a_j)\equiv\alpha$, where $j<h$. Apply Claim 1 to
$z=\varphi_i(a_j)$ and exponent $R_j-1$. There are $0\leq t_j<M$ and
$a_{j+1}\in\mathcal{C}_i$ such that
$\varphi_i^{t_j+1}(a_j)\equiv a_{j+1}$. Put $s_j=t_j+1$ and
$R_{j+1}=R_j-s_j$. Then $1\leq s_j\leq M$,
$R_{j+1}\geq(h-j-1)M+1$, and
$\varphi_i^{R_{j+1}}(a_{j+1})\equiv\alpha$. By (b),
$a_{j+1}=\varphi_i^{s_j}(a_j)$ over $L$, so $a_{j+1}$ is wandering.

This constructs $h+1$ points $a_0,\ldots,a_h$ in $\mathcal{W}_i$.
Hence $a_u=a_v$ for some $u<v$, and the intervening equalities give
$\varphi_i^{s_u+\cdots+s_{v-1}}(a_u)=a_u$, contradicting that $a_u$
is wandering. Thus \((h+1)M\) is a valid bound for the replacement system. By the reduction above, \(B(h+1)M\) is a valid bound for the original system, proving Theorem \ref{thm:intermediate}. \(\square\)

\subsection{Proof of Theorem \ref{thm:split-arithmetic}}

\begin{proof}
The point \(x\), the equations defining \(V\), and all coefficients involved are contained in a common finite extension of \(K\). Replacing \(K\) by that extension preserves the hypotheses and the return set, so we may assume that \(x\) and \(V\) are defined over \(K\).

We first remove preperiodic coordinates by induction on \(g\). For \(g=1\), a preperiodic orbit has an eventually periodic return set. Now suppose, after reindexing, that \(x_g\) is preperiodic, with tail length \(u\) and eventual period \(q\). For \(0\leq r<q\), its last coordinate is constant at the times \(n=u+r+qk\). Intersect \(V\) with the corresponding fiber of the last projection and apply the induction hypothesis to the first \(g-1\) coordinates, the maps \(\varphi_1^q,\ldots,\varphi_{g-1}^q\), and the starting point \(\bigl(\varphi_1^{u+r}(x_1),\ldots,\varphi_{g-1}^{u+r}(x_{g-1})\bigr)\).
An arithmetic progression in \(k\) remains one after the substitution \(n=u+r+qk\), and the finitely many times \(n<u\) are singleton progressions. Thus we may assume that every \(x_i\) is wandering; in particular, the orbit of \(x\) is injective.

Apply Theorem \ref{thm:intermediate} to \(\mathcal T_i=\{x_i\}\). It gives an integer \(M\) and a positive-density set of primes with
\[
 \varphi_i^n(x_i)\not\equiv c\pmod{\mathfrak p}
 \quad\text{for every }c\in\Crit(\varphi_i)\text{ and }n\geq M.
\]
Delete the finitely many primes where some \(\varphi_i\) has bad reduction, where the residue characteristic is at most \(\deg(\varphi_i)\), or where the reduced critical support differs from the reduction of \(\Crit(\varphi_i)\). Fix a remaining prime, choose an extension to \(\overline K\) and the resulting embedding into \(\mathbb C_p\), and view all data there. Put \(y=\Phi^M(x)\).

The orbit of \(y_i\) avoids the residue disk of every attracting
periodic point. Indeed, let \(z_0,\ldots,z_{r-1}\) be an attracting
cycle, and let \(a_j\) denote the derivative of \(\varphi_i\) at \(z_j\),
computed in integral local coordinates. Good reduction gives
\(|a_j|\leq1\) for every \(j\), while attraction gives
\(|a_0\cdots a_{r-1}|<1\). Hence \(\overline{a_j}=0\) for some \(j\),
so \(\overline z_j\in\operatorname{Crit}(\overline{\varphi_i})\). By our
choice of \(\mathfrak p\), there is some
\(c\in\operatorname{Crit}(\varphi_i)\) with
\(\overline z_j=\overline c\). If an iterate of \(y_i\) entered the
residue disk of the cycle, advancing it to the disk of \(z_j\) would
give \(\varphi_i^{M+n}(x_i)\equiv c\pmod{\mathfrak p}\) for some
\(n\geq0\), contradicting the displayed avoidance.

Theorem \ref{thm:padicattracting} now expresses \(V\cap\mathcal O_\Phi(y)\) as a finite union of suborbits. Since the orbit is injective, the corresponding return times form a finite union of arithmetic progressions. Shifting them by \(M\) and adding the first \(M\) times as singleton progressions proves Theorem \ref{thm:split-arithmetic}.
\end{proof}

\section{Relatively \'{E}tale Systems over \texorpdfstring{$\overline{\mathbb Q}$}{Q-bar}}
\label{sec:relative-systems}

We now prove the relative theorem for systems defined over a number field, equivalently over \(\overline{\mathbb Q}\). The relative graph introduced in the Introduction isolates the new hypothesis precisely: the total map is allowed to inherit all ramification of the split base, but has no additional ramification in the relative directions. We first record the formal properties needed to pass to iterates and invariant fibers.

\begin{lemma}[Passage to a positive iterate]\label{lem:iteration}
Let \(F:X\to X\) be an endomorphism and let \(q\geq1\). If the dynamical Mordell--Lang conjecture holds for \(F^q\), then it holds for \(F\).
\end{lemma}

\begin{proof}
Fix \(x\in X\) and a closed subvariety \(V\subseteq X\). For \(0\leq r<q\), put \(S_r=\{m\in\mathbb N_0:(F^q)^m(F^r(x))\in V\}\). By hypothesis, every \(S_r\) is a finite union of arithmetic progressions. Every \(n\in\mathbb N\) has a unique expression \(n=r+qm\) with \(0\leq r<q\), and \(F^{r+qm}(x)=(F^q)^m(F^r(x))\). Hence the return set for \(F\) is \(\bigcup_{r=0}^{q-1}(r+qS_r)\), a finite union of arithmetic progressions.
\end{proof}

\begin{lemma}\label{lem:relative-etale-formal}
Let \(p:X\to B\), \(F:X\to X\), and \(\beta:B\to B\) satisfy \(p\circ F=\beta\circ p\), and suppose that the relative graph \(\delta_F\) is \'{e}tale.

For every \(n\geq1\), the relative graph \(\delta_{F^n}:X\to B\times_{\beta^n,B,p}X\) is \'{e}tale. Suppose additionally that \(\lambda:S\to B\) and \(h:S\to S\) satisfy \(\lambda\circ h=\beta\circ\lambda\). Put \(X_S=S\times_BX\), let \(p_S:X_S\to S\) be the first projection, and define \(F_S(s,x)=(h(s),F(x))\). Then the relative graph of \(F_S\) over \(h\) is \'{e}tale. Finally, if \(\beta\) is \'{e}tale at \(p(x)\) for some \(x\in X\), then \(F\) is \'{e}tale at \(x\).
\end{lemma}

\begin{proof}
For the iteration assertion, put \(Y_n=B\times_{\beta^n,B,p}X\) and write \(\delta_n=\delta_{F^n}\). Thus \(Y_n\) represents pairs \((b,y)\) satisfying \(\beta^n(b)=p(y)\), and \(\delta_n(x)=(p(x),F^n(x))\). Define \(\rho_n:Y_{n+1}\to Y_n\) by \(\rho_n(b,y)=(\beta(b),y)\). This is well-defined because \(\beta^n(\beta(b))=\beta^{n+1}(b)=p(y)\).

Consider \(Y_{n+1}\times_{Y_n,\rho_n,\delta_n}X\). Functorially, a point of this fiber product is a triple \(((b,y),x)\) for which \(\beta(b)=p(x)\) and \(y=F^n(x)\). Conversely, those two identities imply \(\beta^{n+1}(b)=\beta^n(p(x))=p(F^n(x))=p(y)\), so they already ensure that \((b,y)\) belongs to \(Y_{n+1}\). Projection to \((b,x)\) therefore gives a scheme-theoretic identification of this fiber product with \(Y_1\). Under it, the projection to \(Y_{n+1}\) is the map \(\kappa_n:Y_1\to Y_{n+1}\) given by \(\kappa_n(b,x)=(b,F^n(x))\). Hence \(\kappa_n\) is the base change of \(\delta_n\) along \(\rho_n\), and is \'{e}tale whenever \(\delta_n\) is \'{e}tale. Moreover, \(\delta_{n+1}=\kappa_n\circ\delta_1\). Starting from the assumed \'{e}taleness of \(\delta_1=\delta_F\), stability under base change and composition proves the assertion by induction.

The identity \(\lambda\circ h=\beta\circ\lambda\) first shows that \(F_S\) is well-defined: if \(p(x)=\lambda(s)\), then \(p(F(x))=\beta(p(x))=\lambda(h(s))\). Let \(q_1:Y_1\to B\) be the first projection. The target \(S\times_{h,S,p_S}X_S\) of the relative graph of \(F_S\) is canonically isomorphic to \(S\times_{\lambda,B,q_1}Y_1\); explicitly, \((s,(h(s),y))\) corresponds to \((s,(\lambda(s),y))\). Indeed, the condition on either side is \(p(y)=\lambda(h(s))=\beta(\lambda(s))\). Under this identification, the relative graph sends \((s,x)\) to \((s,(\lambda(s),F(x)))\), which is exactly the base change of \(\delta_F:X\to Y_1\) along \(\lambda:S\to B\). It is therefore \'{e}tale.

For the last assertion, let \(\pi_2:Y_1\to X\) be the second projection. The square defining \(Y_1\) identifies \(\pi_2\) as the base change of \(\beta:B\to B\) along \(p:X\to B\). At the point \(\delta_F(x)=(p(x),F(x))\), the corresponding point in the source of \(\beta\) is \(p(x)\). Thus \(\pi_2\) is \'{e}tale at \(\delta_F(x)\) whenever \(\beta\) is \'{e}tale at \(p(x)\). Since \(\delta_F\) is \'{e}tale at \(x\) and \(F=\pi_2\circ\delta_F\), stability under composition gives the conclusion.
\end{proof}

The next lemma is the orbitwise form of the standard \(p\)-adic interpolation argument. Its local input is Poonen's interpolation theorem \cite[Theorem~1]{PoonenInterpolation}.

\begin{lemma}[Interpolation along an \'{e}tale reduced orbit]
\label{lem:orbitwise-interpolation}
Let \(\ell\geq3\) be a prime, let \(\mathcal X\) be a smooth separated scheme of finite type over \(\mathbb Z_\ell\), let \(F:\mathcal X\to\mathcal X\) be a \(\mathbb Z_\ell\)-morphism, and let \(x\in\mathcal X(\mathbb Z_\ell)\). Suppose that \(\overline F\) is \'{e}tale at \(\overline{F^n(x)}\) for every \(n\geq0\). Then, for every closed subvariety \(V\subseteq\mathcal X_{\mathbb Q_\ell}\), the set \(\{n\geq0:F^n(x)\in V\}\) is a finite union of arithmetic progressions.
\end{lemma}

\begin{proof}
The main idea of the proof of the lemma is that the reduced orbit is eventually periodic. On each residue-class tail, we choose analytic coordinates in which a suitable iterate is congruent to the identity modulo \(\ell\). Poonen's theorem then interpolates the orbit, and Strassmann's theorem gives the required zero-set dichotomy.
Write \(x_n=F^n(x)\) and \(\xi_n=\overline{x_n}\). Since \(\mathcal X\) is of finite type over \(\mathbb Z_\ell\), the set \(\mathcal X(\mathbb F_\ell)\) is finite. Reduction commutes with \(F\), so \(\xi_{n+1}=\overline F(\xi_n)\); the sequence is therefore an orbit in a finite set and is eventually periodic. Thus there are \(c\geq0\) and \(q\geq1\) such that \(\xi_{n+q}=\xi_n\) for every \(n\geq c\).

Put \(H=F^q\), and, for \(0\leq j<q\), put \(a_j=x_{c+j}\) and \(\zeta_j=\xi_{c+j}\). Then \(\overline H(\zeta_j)=\zeta_j\). By the chain rule, \(d\overline H_{\zeta_j}\) is the composite of the differentials of \(\overline F\) at \(\xi_{c+j},\ldots,\xi_{c+j+q-1}\). Each factor is an isomorphism by hypothesis, and hence so is \(d\overline H_{\zeta_j}\).

Fix \(j\), and let \(D_j\) be the residue disk consisting of the points of \(\mathcal X(\mathbb Z_\ell)\) reducing to \(\zeta_j\). The equality \(\overline H(\zeta_j)=\zeta_j\) gives \(H(D_j)\subseteq D_j\). Let \(d\) be the relative dimension of \(\mathcal X\) at \(\zeta_j\). If \(d=0\), smoothness makes \(\mathcal X\to\operatorname{Spec}\mathbb Z_\ell\) \'{e}tale on a neighborhood of \(\zeta_j\). Every point reducing to \(\zeta_j\) factors through that neighborhood, and Hensel's lemma gives a unique such lift. Hence \(D_j=\{a_j\}\). It follows from \(H(D_j)\subseteq D_j\) that \(H(a_j)=a_j\), so the progression \(c+j+q\mathbb N\) either consists entirely of return times or contains none.

Assume \(d>0\). Smoothness provides an affine neighborhood \(\mathcal U_j\) of \(\zeta_j\) and an \'{e}tale morphism \(t=(t_1,\ldots,t_d):\mathcal U_j\to\mathbb A^d_{\mathbb Z_\ell}\). The point \(a_j\), and indeed every point of \(D_j\), factors through \(\mathcal U_j\): the inverse image of this open set in \(\operatorname{Spec}\mathbb Z_\ell\) contains the closed point and hence is the whole spectrum. After translating the target by the integral point \(t(a_j)\), we may assume that \(t(a_j)=0\). For each \(u\in\ell\mathbb Z_\ell^d\), Hensel's lemma for the \'{e}tale map \(t\) gives a unique point reducing to \(\zeta_j\) and mapping to \(u\). Consequently, \(t:D_j\to\ell\mathbb Z_\ell^d\) is a bijection. Let \(\iota_j:\mathbb Z_\ell^d\to D_j\) be the inverse parametrization determined by \(t(\iota_j(T))=\ell T\); in particular, \(\iota_j(0)=a_j\).

This parametrization is analytic. The \'{e}tale coordinates identify the completed local ring \(\widehat{\mathcal O}_{\mathcal X,\zeta_j}\) with \(\mathbb Z_\ell[[Z_1,\ldots,Z_d]]\). Thus a regular function on \(\mathcal U_j\) has an expansion \(\sum_\nu a_\nu Z^\nu\) there, with \(a_\nu\in\mathbb Z_\ell\). Pullback by \(\iota_j\) substitutes \(Z=\ell T\); the coefficient of \(T^\nu\) is then \(a_\nu\ell^{|\nu|}\), which tends to zero \(\ell\)-adically as \(|\nu|\) tends to infinity. The pullback therefore belongs to the Tate algebra \(\mathbb Z_\ell\langle T_1,\ldots,T_d\rangle\).

The inclusion \(H(D_j)\subseteq D_j\) means that \(H\) induces an endomorphism of the formal residue disk, and hence defines an analytic self-map \(G_j\) of \(\mathbb Z_\ell^d\) by \(G_j(T)=\ell^{-1}t(H(\iota_j(T)))\). To compute its reduction, write the \(k\)-th coordinate series of this formal endomorphism in the unscaled coordinates, for \(1\leq k\leq d\), as \(\Phi_k(Z)=b_k+\sum_{|\nu|\geq1}a_{k,\nu}Z^\nu\). The coefficients are integral, and \(b_k\in\ell\mathbb Z_\ell\) because \(H(a_j)\in D_j\). After substituting \(Z=\ell T\) and dividing by \(\ell\), all terms of total degree at least two vanish modulo \(\ell\). Hence \(G_j\) reduces to an affine transformation \(T\mapsto A_jT+e_j\). The matrix \(A_j\) represents \(d\overline H_{\zeta_j}\) in the chosen coordinates, so \(A_j\in\operatorname{GL}_d(\mathbb F_\ell)\).

The reduction of \(G_j\) is therefore an element of the finite group \(\operatorname{AGL}_d(\mathbb F_\ell)\). Choose \(s_j\geq1\) divisible by its order, so that \(G_j^{s_j}(T)\equiv T\pmod\ell\) as a congruence of restricted power series. Put \(\Psi_j=G_j^{s_j}\). Since \(\ell\geq3\), the congruence exponent \(1\) satisfies \(1>1/(\ell-1)\), which is the numerical hypothesis in Poonen's interpolation theorem. Applying that theorem to \(\Psi_j\), with initial point \(G_j^t(0)\), gives, for every \(0\leq t<s_j\), a tuple \(\gamma_{j,t}\in\mathbb Z_\ell\langle w\rangle^d\) such that \(\gamma_{j,t}(n)=\Psi_j^n(G_j^t(0))=G_j^{t+s_jn}(0)\) for every \(n\geq0\). Since \(\iota_j\circ G_j=H\circ\iota_j\), we obtain
\[
 F^{c+j+q(t+s_jn)}(x)=\iota_j(\gamma_{j,t}(n))
 \quad\text{for every }n\geq0.
\]

Because \(\mathcal U_j\) is affine, the ideal of \(V\cap(\mathcal U_j)_{\mathbb Q_\ell}\) is generated by finitely many elements of \(\mathcal O(\mathcal U_j)\otimes_{\mathbb Z_\ell}\mathbb Q_\ell\). Choose generators \(f_1,\ldots,f_N\), and multiply each by a power of \(\ell\), when necessary, to make it regular on \(\mathcal U_j\); this does not change their common zero locus on the generic fiber. 

For \(1\leq m\leq N\), the completed-local-ring calculation above gives \(f_m\circ\iota_j\in\mathbb Q_\ell\langle T_1,\ldots,T_d\rangle\). Every coordinate of \(\gamma_{j,t}\) belongs to \(\mathbb Z_\ell\langle w\rangle\) and has Gauss norm at most one, so substitution converges in the Tate algebra. Thus the composites \(h_{m,j,t}=f_m\circ\iota_j\circ\gamma_{j,t}\) belong to \(\mathbb Q_\ell\langle w\rangle\). 

If all of them vanish identically, the entire progression \(c+j+qt+qs_j\mathbb N\) consists of return times. Otherwise, one of them is nonzero, and Strassmann's theorem \cite[Chapter~VI, Section~1.2]{RobertPadicAnalysis} shows that it has only finitely many zeros in \(\mathbb Z_\ell\). The progression then contains only finitely many return times, which may be recorded as singleton progressions.

Finally, every integer at least \(c\) has a unique expression \(c+j+qm\) with \(0\leq j<q\), and every \(m\geq0\) has a unique expression \(t+s_jn\) with \(0\leq t<s_j\). The progressions just considered therefore cover the whole orbit tail. Taking their finite union and adjoining the first \(c\) times proves the lemma.
\end{proof}

\subsection{Proof of the Relative Algebraic Case: Theorem \ref{thm:relative-algebraic}}

\begin{proof}
An arbitrary point of \(X(\overline K)\), together with finitely many equations defining a closed subvariety of \(X_{\overline K}\), is defined over a common finite extension of \(K\). Replacing \(K\) by that extension preserves smoothness, the equivariance identity, and the \'{e}taleness of the relative graph. It is therefore enough to consider \(x\in X(K)\) and a closed subvariety \(V\subseteq X\) defined over \(K\). Write \(p(x)=(a_1,\ldots,a_g)\), and put \(S_F(x,V)=\{n\geq0:F^n(x)\in V\}\). We argue by induction on \(g\). For this induction we allow a smooth quasi-projective scheme with finitely many components and an arbitrary closed subset: the proof below does not use irreducibility, and a closed subset is the finite union of its irreducible components.

Suppose first that some \(a_i\) is preperiodic. After reindexing, take \(i=g\), and choose \(u\geq0\) and \(q\geq1\) such that \(\varphi_g^{n+q}(a_g)=\varphi_g^n(a_g)\) whenever \(n\geq u\). For \(0\leq r<q\), put \(c_r=\varphi_g^{u+r}(a_g)\), \(B_r=(\mathbb P^1)^{g-1}\times\{c_r\}\), \(X_r=p^{-1}(B_r)\), and \(x_r=F^{u+r}(x)\). The point \(c_r\) is fixed by \(\varphi_g^q\), so \(\beta^q\) preserves \(B_r\). The identity \(p\circ F^q=\beta^q\circ p\) then shows that \(F^q\) preserves \(X_r\).

Under the natural identification \(B_r\simeq(\mathbb P^1)^{g-1}\), the restriction of \(\beta^q\) is \(h_r=(\varphi_1^q,\ldots,\varphi_{g-1}^q)\), whose coordinate maps still have degree at least two. If \(\lambda_r:B_r\hookrightarrow B\) is the inclusion, then \(\lambda_r\circ h_r=\beta^q\circ\lambda_r\). Lemma \ref{lem:relative-etale-formal} first makes \(F^q\) relatively \'{e}tale over \(\beta^q\), and then makes \(G_r=F^q|_{X_r}\) relatively \'{e}tale over \(h_r\) by this equivariant base change. The morphism \(X_r\to B_r\) is smooth because it is the base change of \(p\).

For completeness, disconnected fibers cause no problem in the induction. After one further finite extension of \(K\), which neither changes the return set nor the hypotheses, all geometric irreducible components of the finitely many schemes \(X_r\) are defined over \(K\). Smoothness over \(K\) makes these components disjoint and open-and-closed. Let \(C_k\) be the component containing \(G_r^k(x_r)\). The induced self-map of this finite component set makes the sequence \((C_k)\) eventually periodic, so there are \(t\geq0\) and \(s\geq1\) such that \(C_{k+s}=C_k\) for every \(k\geq t\). For \(0\leq h<s\), the map \(G_r^s\) preserves \(C_{t+h}\) and lies over \(h_r^s=(\varphi_1^{qs},\ldots,\varphi_{g-1}^{qs})\). Its relative graph is \'{e}tale: Lemma \ref{lem:relative-etale-formal} gives the relative graph of \(G_r^s\) on \(X_r\), and restricting it to the relevant open-and-closed source and target pieces preserves \'{e}taleness. The induction hypothesis applies to \(G_r^s|_{C_{t+h}}\), the point \(G_r^{t+h}(x_r)\), and the closed subset \(V\cap C_{t+h}\). It follows that each residue-class tail \(t+h+s\mathbb N\) contributes a finite union of arithmetic progressions. Adding the first \(t\) iterates proves the same assertion for the return set of \(G_r\) from \(x_r\).

Every \(n\geq u\) has a unique expression \(n=u+r+qk\), with \(0\leq r<q\) and \(k\geq0\), and then \(F^n(x)=G_r^k(x_r)\). Consequently,
\begin{equation}
 S_F(x,V)=\bigl(S_F(x,V)\cap\{0,\ldots,u-1\}\bigr)
 \cup\bigcup_{r=0}^{q-1}
 \bigl(u+r+qS_{F^q|_{X_r}}(x_r,V\cap X_r)\bigr).
 \label{eq:relative-return-decomposition}
\end{equation}
Each return set on the right is a finite union of arithmetic progressions by induction, and the substitution \(k\mapsto u+r+qk\) carries such a progression to another one. The first set is finite and may be written as a union of singleton progressions. This proves the result in the preperiodic-coordinate case.

We may now suppose that every \(a_i\) is wandering; for \(g=0\) this condition is vacuous, and the following argument supplies the base of the induction. Enlarge \(K\) so that all geometric critical points of the \(\varphi_i\) are rational. This finite extension changes neither wandering nor the return set.

After enlarging a finite set \(S\) of places, we spread all the data
over \(\mathfrak o_{K,S}\).
Put
\(\mathcal B=(\mathbb P^1_{\mathfrak o_{K,S}})^g\), and choose a
quasi-projective model \(\mathcal X\), morphisms
\(\varpi:\mathcal X\to\mathcal B\) and
\(\mathcal F:\mathcal X\to\mathcal X\), extensions of the maps
\(\varphi_i\), a section extending \(x\), and the schematic closure
\(\mathcal V\) of \(V\).
Writing \(\mathcal\beta\) for the product of the extended coordinate
maps, we retain the identity
\(\varpi\circ\mathcal F=\mathcal\beta\circ\varpi\).
Smoothness and \'{e}taleness are open conditions for morphisms of
finite presentation.
Since \(\varpi\) is smooth and the relative graph of \(\mathcal F\)
over \(\mathcal\beta\) is \'{e}tale on the generic fiber, enlarging
\(S\) once more makes both properties hold over
\(\mathfrak o_{K,S}\).
The scheme \(\mathcal X\) is then smooth over
\(\mathfrak o_{K,S}\), because both \(\varpi\) and
\(\mathcal B\to\operatorname{Spec}\mathfrak o_{K,S}\) are smooth;
it is also separated and of finite type because it is
quasi-projective.

For each \(i\), choose a homogeneous presentation of \(\varphi_i\).
After inverting its resultant, this presentation extends to every
remaining fiber without a drop in degree.
Its homogeneous Wronskian defines the ramification divisor.
Since all geometric critical points are \(K\)-rational, factor the
Wronskian over \(K\) into the distinct linear forms defining those
points, with their ramification multiplicities.
Choose primitive integral coordinates for the linear forms and enlarge
\(S\) so that the scalar in the factorization is a unit, the
determinant of every pair of distinct linear forms is a unit, and no
remaining residue characteristic divides any \(\deg(\varphi_i)\).
The reduced linear forms are then distinct and the Wronskian
factorization remains nonzero after reduction.
Consequently, at every remaining prime, each \(\varphi_i\) has good
reduction of unchanged degree, and the support of its reduced
ramification divisor is exactly the reduction of its geometric
critical set.

If \(g>0\), apply Theorem \ref{thm:intermediate} with \(\mathcal T_i=\{a_i\}\). It gives \(M\geq1\) and a positive-density set \(\mathcal P\) of places such that \(\varphi_i^{M+n}(a_i)\) is not congruent to any critical point of \(\varphi_i\), for every \(i\) and every \(n\geq0\); here we have written the theorem's index \(m\geq M\) as \(m=M+n\). The proof of Theorem \ref{thm:intermediate} obtains this set from Chebotarev sets by finite removal, so \(\mathcal P\) has positive Dirichlet density. If \(g=0\), put \(M=0\) and let \(\mathcal P\) be the set of all finite places. Remove from \(\mathcal P\) the finite set \(S\), the ramified places, the places above \(2\), and, when \(g>0\), the places whose residue characteristic is at most \(\max_i\deg(\varphi_i)\). Only finitely many places have been removed, so \(\mathcal P\) still has positive Dirichlet density.

The places of \(K\) having absolute residue degree greater than one have Dirichlet density zero. Indeed, if such a place lies over \(\ell\), its norm is at least \(\ell^2\), and the contribution of all such places to the prime Dirichlet series is bounded near \(s=1\) by \([K:\mathbb Q]\sum_\ell\ell^{-2s}\), which remains bounded as \(s\to1^+\). We may therefore choose \(v\in\mathcal P\) of absolute residue degree one. By the preceding exclusions, \(v\) is unramified over a rational prime \(\ell\geq3\). Its ramification index and residue degree are both one, so \([K_v:\mathbb Q_\ell]=1\); hence \(K_v=\mathbb Q_\ell\) and \(\mathfrak o_{K_v}=\mathbb Z_\ell\).

Base-change the integral model to \(\mathbb Z_\ell\), and put \(z=F^M(x)\). This is an integral point because both \(x\) and \(F\) extend to the model. For every \(n\geq0\), equivariance gives \(p(F^n(z))=\beta^{M+n}(p(x))\). The \(i\)-th coordinate of its reduction is therefore \(\overline{\varphi_i^{M+n}(a_i)}\), which is not the reduction of any critical point. By the construction of the ramification model, \(\overline{\varphi_i}\) is \'{e}tale at that coordinate. Thus the product map \(\overline\beta\) is \'{e}tale at every point of the reduced orbit of \(p(z)\); when \(g=0\), this is automatic because \(\beta\) is the identity of a point.

The special fiber of the relative graph of \(F\) remains \'{e}tale by base change. If \(z_n=\overline{F^n(z)}\), then \(\overline p(z_n)=\overline\beta^{\,n}(\overline p(z))\). Applying the last assertion of Lemma \ref{lem:relative-etale-formal} on the special fiber, at each \(z_n\), shows that \(\overline F\) is \'{e}tale at every point of the reduced orbit of \(z\). The base-changed \(\mathcal X\) is smooth, separated, and of finite type over \(\mathbb Z_\ell\), so Lemma \ref{lem:orbitwise-interpolation}, applied to the base change of \(V\), proves that the \(p\)-adic return set \(\{n\geq0:F^n(z)\in V_{\mathbb Q_\ell}\}\) is a finite union of arithmetic progressions. To compare this with the original return set, fix \(n\) and choose an affine \(K\)-neighborhood \(W\) of \(F^{M+n}(x)\), with generators \(f_1,\ldots,f_r\) for the ideal of \(V\cap W\). The point belongs to \(V\) exactly when every \(f_j(F^{M+n}(x))\) is zero in \(K\). Since \(K\hookrightarrow K_v=\mathbb Q_\ell\) is injective, these equalities hold exactly when their images in \(K_v\) vanish. The \(p\)-adic return set is therefore exactly the return set of the original orbit after time \(M\). Translating these progressions by \(M\) and adjoining the return times smaller than \(M\) as singleton progressions proves the theorem.
\end{proof}

\part{Simultaneous Specialization, prescribed-reduction embeddings, and Results over Complex Numbers}

\section{Specialization and the Passage to \texorpdfstring{$\mathbb C$}{C}}
\label{sec:complex-transfer}

The preceding theorem applies to systems defined over \(\overline{\mathbb Q}\), whereas a complex system need not descend to \(\overline{\mathbb Q}\). We first specialize the wandering base coordinates to algebraic data in order to select a favorable prime. We then lift the same reduction through an injective \(p\)-adic realization of the original coefficient field. Only the second, injective map is used to compare return sets.

\begin{lemma}[Simultaneous specialization of wandering points]
\label{lem:simultaneous-specialization}
Let \(k\) be a number field, let \(T\) be a geometrically integral \(k\)-variety with function field \(K=k(T)\), and let \(U\subseteq T\) be nonempty and open. For \(1\leq i\leq r\), let \(\varphi_i:\mathbb P^1_K\to\mathbb P^1_K\) have degree at least \(2\), and let \(\mathcal T_i\subseteq\mathbb P^1(K)\) be finite and contain no \(\varphi_i\)-preperiodic point. Suppose that the maps and marked points extend over \(U\), without a drop in degree. Then there are a finite extension \(L/k\) and \(t_0\in U(L)\) such that every specialization \(a_{t_0}\), with \(a\in\mathcal T_i\), is wandering under \(\varphi_{i,t_0}\). Any fixed finite collection of nonvanishing conditions on the spread may be retained by shrinking \(U\).
\end{lemma}

\begin{proof}
Fix an algebraic closure of \(k\) inside \(\overline{\mathbb Q}\). Geometric integrality of \(T\) implies that \(K\overline{\mathbb Q}\) is a field, so the base-change map from \(K\) is injective. Preperiodicity is preserved and reflected by this extension: it is the existence of two indices \(m>n\) for which the equality \(\varphi_i^m(a)=\varphi_i^n(a)\) holds. We may therefore extend the constant field from \(k\) to \(\overline{\mathbb Q}\) without changing the hypothesis on any of the marked points. At the end, the coordinates of the specialization point involve only finitely many algebraic numbers, and hence lie in a finite extension of \(k\).

We use the transcendence-degree induction of \cite[Lemma~3.3]{GhiocaXieSkewLinear} and \cite[Claim~3.4]{GhiocaXieSkewLinear}, which extends the one-parameter argument of \cite[Proposition~6.2]{GTZ08}. We spell out why the induction can be carried out simultaneously for the present finite family. Its one-step assertion is the following: if the current coefficient field is the function field \(E(C)\) of a smooth projective curve over a field \(E\) arising at the preceding stage, then one can choose a closed point of \(C\), with residue field finite over \(E\), at which every current marked point remains wandering. The maps, points, and any prescribed nonzero functions are first spread over a common nonempty open subset of \(C\), so that specialization is defined and no degree drops there.

Use the global-height-field structure from the cited induction. When \(E\) has positive transcendence degree, choose a normal projective model together with an ample polarization; weighting the divisorial valuations by the resulting intersection degrees gives the product formula required for a global height field. At the bottom stage, the finitely many data are first descended from \(\overline{\mathbb Q}\) to a number field and its ordinary Weil height is used. Fix an ample divisor \(D\) on \(C\), and normalize its associated height by putting \(h_C=h_D/\deg D\). For a pair \((\varphi_i,a)\), put \(d_i=\deg(\varphi_i)\), and let \(\widehat h_{\varphi_i}(a)\) be its geometric canonical height over \(E(C)\). Suppose first that this height is positive. The specialization theorem \cite[Theorem~4.1]{CallSilverman1993}, in the global-height-field form used in \cite[Claim~3.4]{GhiocaXieSkewLinear}, gives
\[
 \lim_{h_C(t)\to\infty}
 \frac{\widehat h_{\varphi_{i,t}}(a_t)}{h_C(t)}
 =\widehat h_{\varphi_i}(a)>0.
\]
The quotient is consequently positive for every good specialization of sufficiently large height. Since a preperiodic point has canonical height zero, each such \(a_t\) is wandering.

Suppose instead that \(\widehat h_{\varphi_i}(a)=0\). We may extend the constants to an algebraic closure of \(E\): the equalities defining preperiodicity are preserved and reflected by a field extension, and vanishing of canonical height is unchanged. Baker's theorem \cite[Corollary~1.8]{Baker2009} applies over this algebraically closed constant field. Because \(a\) is wandering, it forces \(\varphi_i\) to be isotrivial. The conjugacy and the resulting constant map involve only finitely many coefficients, so they descend to a finite extension of \(E\). After replacing \(C\) by the corresponding finite cover \(C'\), there is a linear-fractional transformation \(\mu\) for which \(\psi=\mu\circ\varphi_i\circ\mu^{-1}\) is defined over that finite extension of the constant field. Put \(b=\mu(a)\). Finite base change multiplies the geometric height only by a positive normalization factor, while conjugacy preserves its vanishing; hence \(\widehat h_\psi(b)=0\). If \(b\) were nonconstant, it would define a nonconstant morphism \(b:C'\to\mathbb P^1\), and \(\deg(\psi^n\circ b)=d_i^n\deg(b)\). The defining limit for the geometric canonical height would then be a positive multiple of \(\deg(b)\), contrary to its vanishing. Thus \(b\) is constant. It is wandering because conjugacy preserves preperiodicity and \(a\) is wandering. Remove the finitely many points of \(C'\) at which the chosen models of \(\mu\), \(\mu^{-1}\), \(\varphi_i\), or \(a\) do not specialize with the identities \(\psi=\mu\varphi_i\mu^{-1}\) and \(b=\mu(a)\) intact. At every remaining point the specialized pair is conjugate to the same constant pair \((\psi,b)\), and hence \(a_t\) is wandering.

There are only finitely many pairs. Take the normalization of \(C\) in the compositum of the finite covering function fields required by the height-zero pairs, and choose a geometrically integral component dominating \(C\), enlarging \(E\) by its finite constant field if necessary. Pullback to this common cover preserves positivity of geometric canonical height. Reapply the specialization theorem on the common cover for the positive-height pairs, express the resulting finitely many bounds using one ample height, and take their maximum. For the height-zero pairs, remove the union of the finitely many bad-specialization sets. The zero and pole loci of all prescribed nonzero functions, as well as the complement of the original open set, are added to the same exceptional set. To see that the remaining points have unbounded height, choose a nonconstant finite morphism \(\pi\) from the common covering curve to \(\mathbb P^1_E\). Functoriality gives \(h_{\pi^*\mathcal O(1)}(t)=h_{\mathbb P^1}(\pi(t))+O(1)\), while the height of the positive-degree divisor \(\pi^*\mathcal O(1)\) is comparable, up to positive scaling and \(O(1)\), with the chosen ample height on the curve. We may therefore lift algebraic points of \(\mathbb P^1\) of arbitrarily large height that avoid the image of the finite exceptional set, and choose one lift beyond all the preceding bounds. Its residue field is finite over \(E\), and every member of the finite family remains wandering there.

Finally, we iterate the one-step construction by induction on \(d=\operatorname{trdeg}_{\overline{\mathbb Q}}(K\overline{\mathbb Q})\), as in \cite[Lemma~3.3]{GhiocaXieSkewLinear}. Before beginning, replace \(U\), if necessary, by a nonempty affine principal open \(U_0=D(f)\subseteq U\). Its coordinate ring is finitely generated; we require at every step that its algebra generators, the coefficients defining all maps and marked points, and the inverse of \(f\) remain regular. This ensures that the successive evaluations compose to a homomorphism from the coordinate ring of \(U_0\), and hence ultimately define a point of the original \(U\).

Suppose \(d>0\), and choose \(u_1,\ldots,u_{d-1}\) as the first members of a separating transcendence basis of the current coefficient field over \(\overline{\mathbb Q}\). Let \(E\) be the relative algebraic closure of \(\overline{\mathbb Q}(u_1,\ldots,u_{d-1})\) there. Then \(E\) is finitely generated of transcendence degree \(d-1\), it is algebraically closed in the current field, and the latter is a one-variable function field over \(E\). It is therefore the function field of a smooth projective geometrically integral curve over \(E\). Apply the simultaneous one-step construction on this curve, outside the poles of the finitely many retained algebra generators and denominators. The resulting evaluation has residue field \(E'\) finite over \(E\), so \(E'\) is again finitely generated over \(\overline{\mathbb Q}\), now of transcendence degree \(d-1\). Because \(\overline{\mathbb Q}\) is algebraically closed, \(E'\) is the function field of a geometrically integral \(\overline{\mathbb Q}\)-variety. The retained nonzero elements, including the homogeneous resultants, define a nonempty open on this variety where the maps have their original degrees and all marked points are defined. Their images remain wandering there. We may therefore apply the induction hypothesis to this specialized finite family and to the image of the coordinate ring. Composing the two evaluations removes all remaining transcendence parameters. The base case \(d=0\) is immediate because a finitely generated algebraic extension of \(\overline{\mathbb Q}\) equals \(\overline{\mathbb Q}\). The composite thus gives a point of \(U_0(\overline{\mathbb Q})\). Its finitely many coordinates generate a finite extension \(L/k\), and it therefore defines the required point \(t_0\in U(L)\).
\end{proof}

The preceding specialization selects a favorable residue point but is
not injective. The next lemma realizes the original coefficient field \(p\)-adically while retaining that prescribed reduction.

\begin{lemma}[Injective lift with prescribed reduction]
\label{lem:prescribed-reduction}
Let \(L\) be a number field, let \(\mathcal U=\operatorname{Spec}\mathcal A\) be a smooth integral affine scheme of finite type over \(\mathfrak o_{L,S}\), and assume that its generic fiber \(U\) is geometrically integral. Let \(t_0\in\mathcal U(\mathfrak o_{L,S})\), and let \(v\notin S\) be unramified over a rational prime \(\ell\) and have absolute residue degree one. Then there is an injective homomorphism \(\iota:\mathcal A\to\mathbb Z_\ell\) whose restriction to \(\mathfrak o_{L,S}\) is induced by \(v\) and whose reduction modulo \(\ell\) equals the reduction of \(t_0\). Consequently, \(\iota\) extends to an embedding \(L(U)\hookrightarrow\mathbb Q_\ell\).
\end{lemma}

\begin{proof}
Since \(v\) is unramified over \(\ell\) and has absolute residue degree one, its ramification index and residue degree are both one. Thus \([L_v:\mathbb Q_\ell]=1\), so \(L_v=\mathbb Q_\ell\) and \(\mathfrak o_{L_v}=\mathbb Z_\ell\). We use the corresponding embedding \(\mathfrak o_{L,S}\hookrightarrow\mathbb Z_\ell\), and put \(\mathcal B=\mathcal A\otimes_{\mathfrak o_{L,S}}\mathbb Z_\ell\).

The ring \(\mathcal B\) is a domain. Smoothness makes it flat, and hence torsion-free, over \(\mathbb Z_\ell\), so it injects into \(\mathcal B[1/\ell]\). The latter is the coordinate ring of the base change of \(U\) to \(\mathbb Q_\ell\), which is integral because \(U\) is geometrically integral. A zero divisor in \(\mathcal B\) would remain one after inverting \(\ell\), unless one of its factors became zero; torsion-freeness rules out the latter possibility. Hence \(\mathcal B\) is integral.

Let \(\overline t_0\) be the reduction of the \(\mathbb Z_\ell\)-point obtained from \(t_0\), let \(\mathcal D\) be the set of points of \(\mathcal U(\mathbb Z_\ell)\) reducing to \(\overline t_0\), and put \(d=\dim U\). If \(d=0\), the coordinate ring of the integral affine variety \(U\) is a finite field extension of \(L\). Geometric integrality forces this extension to be \(L\) itself, so \(U=\operatorname{Spec}L\). The section \(t_0:\mathcal A\to\mathfrak o_{L,S}\) is injective because it becomes the identity after passing to the generic fiber. Composing it with \(\mathfrak o_{L,S}\hookrightarrow\mathbb Z_\ell\) gives the required map.

Assume \(d>0\). By smoothness, there are an affine neighborhood \(\mathcal W\) of \(\overline t_0\) and an \'{e}tale morphism \(z=(z_1,\ldots,z_d):\mathcal W\to\mathbb A^d_{\mathbb Z_\ell}\). Every \(\mathbb Z_\ell\)-point reducing to \(\overline t_0\) factors through \(\mathcal W\), since the inverse image of \(\mathcal W\) in the local scheme \(\operatorname{Spec}\mathbb Z_\ell\) contains the closed point. After translating the target by \(z(t_0)\), we may assume that \(z(t_0)=0\). For every \(u\in\ell\mathbb Z_\ell^d\), Hensel's lemma for the \'{e}tale map \(z\) gives a unique point reducing to \(\overline t_0\) and mapping to \(u\). Thus \(z\) identifies \(\mathcal D\) with \((\ell\mathbb Z_\ell)^d\). With the transported \(\ell\)-adic metric, \(\mathcal D\) is therefore a nonempty complete metric space.

The ring \(\mathcal A\) is countable because it is finitely generated over the countable ring \(\mathfrak o_{L,S}\). For every nonzero \(f\in\mathcal A\), put \(Z_f=\{t\in\mathcal D:f(t)=0\}\). Continuity of evaluation in the \'{e}tale coordinates makes \(Z_f\) closed. We claim that it has empty interior.

Let \(\eta:\mathbb Z_\ell^d\to\mathcal D\) be the parametrization characterized by \(z(\eta(T))=\ell T\). The \'{e}tale coordinates identify the completed local ring at \(\overline t_0\) with \(\mathbb Z_\ell[[Z_1,\ldots,Z_d]]\). If the image of \(f\) there is \(\sum_\nu a_\nu Z^\nu\), then \(f(\eta(T))=\sum_\nu a_\nu\ell^{|\nu|}T^\nu\), an element of \(\mathbb Z_\ell\langle T_1,\ldots,T_d\rangle\). If \(Z_f\) contained an open subset of \(\mathcal D\), this analytic function would vanish on an open subpolydisk. Applying the one-variable nonarchimedean identity theorem successively in the variables would make every coefficient \(a_\nu\ell^{|\nu|}\) zero, and hence would make the image of \(f\) zero in the completed local ring.

That is impossible for a nonzero \(f\). To make the injections explicit, put \(R=\mathfrak o_{L,S}\), \(R_v=\mathfrak o_{L,v}\), and \(\mathcal A_v=\mathcal A\otimes_RR_v\). The map \(\mathcal A\to\mathcal A_v\) is injective because it is a localization of the domain \(\mathcal A\). The completion map \(R_v\to\mathbb Z_\ell\) is faithfully flat, so \(\mathcal A_v\to\mathcal A_v\otimes_{R_v}\mathbb Z_\ell=\mathcal B\) is also injective. Next, the domain \(\mathcal B\) injects into its localization \(\mathcal B_{\overline t_0}\), and the Krull intersection theorem injects this Noetherian local domain into its completion. Thus every nonzero element of \(\mathcal A\) remains nonzero in the completed local ring. This proves that each \(Z_f\) is closed and nowhere dense.

The sets \(Z_f\), as \(f\) ranges over the nonzero elements of \(\mathcal A\), are a countable family of closed nowhere-dense subsets of the complete metric space \(\mathcal D\). The Baire category theorem therefore gives \(t_*\in\mathcal D\setminus\bigcup_{0\ne f\in\mathcal A}Z_f\). Evaluation at \(t_*\) defines a homomorphism \(\iota:\mathcal A\to\mathbb Z_\ell\), and the choice of \(t_*\) says precisely that its kernel is zero. Its restriction to \(\mathfrak o_{L,S}\) is the embedding induced by \(v\). Moreover, \(t_*\) and \(t_0\) lie in the same residue disk, so their induced maps \(\mathcal A\to\mathbb F_\ell\) agree.

Finally, \(\iota\) extends by localization to an embedding \(\operatorname{Frac}(\mathcal A)\hookrightarrow\mathbb Q_\ell\). Passing from \(\mathcal A\) to its generic fiber does not change its fraction field, so \(\operatorname{Frac}(\mathcal A)=L(U)\), as required.
\end{proof}

\begin{remark}\label{rem:two-specializations}
The point \(t_0\) in Lemma \ref{lem:simultaneous-specialization} is used only to choose a favorable arithmetic reduction. Its return set is never identified with the return set of the original complex system. Lemma \ref{lem:prescribed-reduction} instead gives an injective embedding of the original coefficient field while preserving the chosen reduction, and injectivity is what permits the final comparison of return sets.
\end{remark}

\begin{proof}[Proof of Theorem \ref{thm:relative-complex}]
Fix \(x\in X(\mathbb C)\) and a closed subvariety \(V\subseteq X\), and write \(p(x)=(a_1,\ldots,a_g)\). We first remove preperiodic base coordinates by induction on \(g\). If \(g=0\), then \(B\) is a point and the fiber product defining the relative graph is canonically \(X\); under this identification \(\delta_F\) is \(F\). Hence \(F\) is \'{e}tale, and the theorem follows from \cite{Bell2010}.

Suppose \(g>0\) and, after reindexing, that \(a_g\) is preperiodic. Choose \(u\geq0\) and \(q\geq1\) such that \(\varphi_g^{n+q}(a_g)=\varphi_g^n(a_g)\) for every \(n\geq u\). For \(0\leq r<q\), put \(c_r=\varphi_g^{u+r}(a_g)\), \(B_r=(\mathbb P^1)^{g-1}\times\{c_r\}\), \(X_r=p^{-1}(B_r)\), and \(x_r=F^{u+r}(x)\). The equality \(\varphi_g^q(c_r)=c_r\) makes \(G_r=F^q|_{X_r}\) a self-map of \(X_r\). Lemma \ref{lem:relative-etale-formal}, applied first to \(F^q\) and then to the equivariant inclusion \(B_r\hookrightarrow B\), makes \(G_r\) relatively \'{e}tale over \(h_r=(\varphi_1^q,\ldots,\varphi_{g-1}^q)\). The map \(X_r\to B_r\) is smooth by base change, so \(X_r\) is regular and its finitely many irreducible components are disjoint and open-and-closed.

Let \(C_k\) be the component containing \(G_r^k(x_r)\). The sequence \((C_k)\) is eventually periodic in this finite set, so choose \(t\geq0\) and \(s\geq1\) with \(C_{k+s}=C_k\) for every \(k\geq t\). For each \(0\leq h<s\), the map \(G_r^s\) preserves \(C_{t+h}\) and lies over \(h_r^s=(\varphi_1^{qs},\ldots,\varphi_{g-1}^{qs})\). Lemma \ref{lem:relative-etale-formal} gives an \'{e}tale relative graph for \(G_r^s\), and its restriction to the relevant open-and-closed source and target pieces is the relative graph of \(G_r^s|_{C_{t+h}}\). The induction hypothesis therefore applies to this component and to every irreducible component of \(V\cap C_{t+h}\). It gives the desired return-set conclusion on each tail \(t+h+s\mathbb N\); the first \(t\) iterates contribute only singleton progressions. Thus the return set for \(G_r\) from \(x_r\) is a finite union of arithmetic progressions. Finally, every \(n\geq u\) has a unique expression \(n=u+r+qk\), and then \(F^n(x)=G_r^k(x_r)\). The decomposition \eqref{eq:relative-return-decomposition} proves the result. We may consequently assume that every \(a_i\) is wandering.

Choose a finitely generated field \(K_0\subseteq\mathbb C\) over which \(B,X,p,F,x,V\), the maps \(\varphi_i\), the relative graph, and all compatibility identities are defined. Let \(k\) be the relative algebraic closure of \(\mathbb Q\) in \(K_0\). Because \(K_0\) is finitely generated over \(\mathbb Q\), the algebraic extension \(k/\mathbb Q\) is finite. The field \(k\) is algebraically closed in \(K_0\), and every extension in characteristic zero is separable, so \(K_0/k\) is regular. If \(K_0=k\), all the data are defined over a number field and Theorem \ref{thm:relative-algebraic} applies after base change to an algebraic closure of \(K_0\). Its conclusion gives the complex return set as well: membership in \(V\) is described by equations over \(K_0\), and the injection \(K_0\hookrightarrow\mathbb C\) preserves and reflects their vanishing. We may therefore assume that \(K_0\) has positive transcendence degree over \(k\).

Regularity gives a geometrically integral \(k\)-variety with function field \(K_0\), and characteristic zero permits us to replace it by a smooth nonempty open \(T\). After a further shrinking to a nonempty affine open \(U\), all the data spread over \(U\). Indeed, schemes, morphisms, sections, equations, and their compatibility identities involve only finitely many coefficients. The relevant smooth loci and the \'{e}tale locus of the relative graph are open and contain the generic fiber. We may therefore shrink \(U\) so that the models of \(B\) and \(X\) are smooth over \(U\), the model of \(p\) is a smooth morphism, and the relative graph is an \'{e}tale morphism. We also invert the homogeneous resultants of the \(\varphi_i\), which prevents a drop in degree on any fiber. The point \(x\), its coordinate sections \(a_i\), and the schematic closure of \(V\) extend over the same \(U\).

Apply Lemma \ref{lem:simultaneous-specialization} to the pairs \((\varphi_i,a_i)\). It gives a finite extension \(L/k\) and \(t_0\in U(L)\) such that every \(a_{i,t_0}\) is wandering under \(\varphi_{i,t_0}\). Enlarge \(L\) so that all critical points of the specialized maps are \(L\)-rational, and choose an embedding \(L\hookrightarrow\mathbb C\) extending the fixed embedding of \(k\). Since \(K_0/k\) is regular, \(K_0\) and \(L\) are linearly disjoint over \(k\). The multiplication map \(K_0\otimes_kL\to\mathbb C\) is therefore injective, and the fraction field of its image is the compositum \(K=K_0L\subseteq\mathbb C\). Geometric integrality of \(U\) makes \(U_L\) integral, and its function field is precisely \(K\).

After enlarging a finite set \(S\) of places of \(L\), choose a smooth integral affine model \(\mathcal U=\operatorname{Spec}\mathcal A\) of \(U_L\) over \(\mathfrak o_{L,S}\) such that \(t_0\) extends to an \(\mathfrak o_{L,S}\)-point. Spread the entire system over \(\mathcal U\). We may arrange that the total space is separated and of finite type, that the model of \(p\) is smooth, that the relative graph is \'{e}tale, and that the maps \(\varphi_i\) have unchanged degree. Since the split base is smooth over \(\mathcal U\), smoothness of \(p\) also makes the total space smooth over \(\mathcal U\). The point \(x\) extends as a section and \(V\) by schematic closure.

At this stage we need a ramification condition only for the maps specialized at \(t_0\), not a splitting of the generic critical locus of the family. Choose homogeneous presentations of the maps \(\varphi_{i,t_0}\). Their ramification divisors are defined by the corresponding homogeneous Wronskians, which split over \(L\) into the linear forms defining the chosen \(L\)-rational critical points, with their multiplicities. Choose primitive integral coordinates for these points. By adding finitely many places to \(S\), we make the resultants and the scalar factors in these factorizations units, keep the determinants of distinct linear factors invertible, and exclude residue characteristics dividing any \(\deg(\varphi_i)\). For every place outside \(S\), the reduced Wronskian is then nonzero and the support of the ramification divisor of \(\overline{\varphi_{i,t_0}}\) is exactly the set of reductions of the characteristic-zero critical points.

Apply Theorem \ref{thm:intermediate} over \(L\) to the specialized maps and to \(\mathcal T_i=\{a_{i,t_0}\}\). We obtain \(M\geq1\) and a positive-density set \(\mathcal P\) of places \(v\) such that \(\varphi_{i,t_0}^{M+n}(a_{i,t_0})\) is not congruent to any critical point of \(\varphi_{i,t_0}\), for every \(i\) and every \(n\geq0\). As in the arithmetic proof, the construction of this set from Chebotarev sets shows that it has positive Dirichlet density. Remove from \(\mathcal P\) the places in \(S\), the ramified places, the places above \(2\), and the places whose residue characteristic is at most \(\max_i\deg(\varphi_i)\). These are finite removals. Places of \(L\) having absolute residue degree greater than one have Dirichlet density zero: their prime Dirichlet series is bounded near \(s=1\) by \([L:\mathbb Q]\sum_\ell\ell^{-2s}\). We may therefore choose a remaining \(v\in\mathcal P\) of absolute residue degree one, above a rational prime \(\ell\geq3\). Since \(v\) is unramified and has residue degree one, \([L_v:\mathbb Q_\ell]=1\), so \(L_v=\mathbb Q_\ell\) and \(\mathfrak o_{L_v}=\mathbb Z_\ell\).

Lemma \ref{lem:prescribed-reduction} gives an injection \(\iota_{\mathcal A}:\mathcal A\to\mathbb Z_\ell\), which extends to a field embedding \(\iota:K=\operatorname{Frac}(\mathcal A)\hookrightarrow\mathbb Q_\ell\). Let \(\tau_0:\mathcal A\to\mathbb F_\ell\) be evaluation at \(t_0\) followed by reduction modulo \(v\). The prescribed-reduction property says that \(\tau_0\) is also the reduction of \(\iota_{\mathcal A}\). Pulling an object in the integral family back through \(t_0\) and then reducing therefore gives the same special fiber as pulling it back through \(\iota_{\mathcal A}\) and reducing modulo \(\ell\). In particular, the two procedures give the same reduced maps and marked points. The maps and marked points extend over \(\mathcal A\), so every iterate is an \(\mathcal A\)-section. Let \(a_{i,\iota}^{(m)}\in\mathbb P^1(\mathbb Z_\ell)\) be the pullback through \(\iota_{\mathcal A}\) of the section \(\varphi_i^m(a_i)\). Since iteration is formed by composition in the integral family, for every \(i\) and \(n\geq0\) the reduction of \(a_{i,\iota}^{(M+n)}\) equals \(\overline{\varphi_{i,t_0}^{M+n}(a_{i,t_0})}\). The latter point avoids the reductions of all critical points of \(\varphi_{i,t_0}\). By the ramification-support condition, it lies outside the ramification locus of the reduced map, which is also the reduction of the map obtained through \(\iota_{\mathcal A}\). Thus the latter map is \'{e}tale at each of the points \((\overline{a_{1,\iota}^{(M+n)}},\ldots,\overline{a_{g,\iota}^{(M+n)}})\), for \(n\geq0\), of the reduced base orbit. Taking the product over \(i\) shows that the reduction of \(\beta\) is \'{e}tale at every one of these points.

Pull the whole integral system back through \(\iota_{\mathcal A}\), and continue to denote its maps by \(p,F\), and \(\beta\). The resulting scheme is smooth, separated, and of finite type over \(\mathbb Z_\ell\), and its relative graph remains \'{e}tale by base change. Because \(x\) is a section of the integral family, \(z=F^M(x)\) belongs to \(\mathcal X(\mathbb Z_\ell)\). Put \(z_0=\overline z\) and \(z_n=\overline F^{\,n}(z_0)=\overline{F^n(z)}\). Then \(\overline p(z_n)=\overline\beta^{\,n}(\overline p(z_0))\), and \(\overline p(z_0)=(\overline{a_{1,\iota}^{(M)}},\ldots,\overline{a_{g,\iota}^{(M)}})\) is the \(n=0\) point of the reduced base orbit just described. The base map is \'{e}tale at this point and at all of its subsequent iterates. Applying the last assertion of Lemma \ref{lem:relative-etale-formal} on the special fiber, at each \(z_n\), shows that \(\overline F\) is \'{e}tale at every point of the reduced orbit of \(z\). Lemma \ref{lem:orbitwise-interpolation}, applied to the generic fiber of the spread of \(V\), therefore shows that \(\{n\geq0:F^n(z)\in V_{\mathbb Q_\ell}\}\) is a finite union of arithmetic progressions.

Put \(w_n=F^{M+n}(x)\in X(K)\).
Both \(K\hookrightarrow\mathbb C\) and
\(\iota:K\hookrightarrow\mathbb Q_\ell\) are field embeddings and
hence faithfully flat.
Therefore \(w_n\in V_K\) if and only if its complex realization lies
in \(V_{\mathbb C}\), and if and only if its \(p\)-adic realization
lies in \(V_{\mathbb Q_\ell}\).
Thus the interpolated return set is exactly the return set of the
original complex orbit after time \(M\), which completes the proof.
\end{proof}

\section{Stable Mixed Models for Products of Curves over \texorpdfstring{$\mathbb C$}{C}}
\label{sec:stable-structure}

We now describe the stable dynamics of an endomorphism of a product of complex curves. An empty product is understood as a point. A stable image need not itself be a product, so the appropriate conclusion is an equivariant surjective parametrization by a product.

The geometric idea is expressed by a dependency graph. We draw an arrow \(i\to j\) when the \(i\)-th output coordinate
depends nontrivially on the \(j\)-th input coordinate. A nonconstant map to a curve of genus at least \(2\) depends on only one curve factor, and a nonconstant map from a product of projective lines to a projective line depends on only one rational factor. Thus the relevant dependency graphs among the hyperbolic coordinates and among the rational coordinates that depend on rational inputs have out-degree at most one; rational outputs factoring through the positive-genus block are treated separately. The resulting directed graph consists of cycles with finite trees feeding into them. After sufficiently many iterations, the tree coordinates are forgotten; only the cycles and the positive-genus stable image remain.

\subsection{Maps from products to curves}

\begin{lemma}\label{lem:product-to-hyperbolic}
Let \(C_1,\ldots,C_n\) and \(D\) be smooth projective connected complex curves with \(g(D)\geq2\). Every nonconstant morphism \(u:\prod_{j=1}^n C_j\to D\) has the form \(u=v\circ\operatorname{pr}_j\) for some \(j\) and some nonconstant morphism \(v:C_j\to D\).
\end{lemma}

\begin{proof}
Since \(u\) is nonconstant in characteristic zero, its differential is nonzero at some point. The tangent space of the product is the direct sum of the tangent spaces of its factors. Consequently, for some index \(j\), the restriction of \(u\) to a fiber obtained by varying only the \(j\)-th coordinate is nonconstant.

Set \(T=\prod_{i\ne j}C_i\), and, for \(t\in T\), let \(u_t:C_j\to D\) be the restriction of \(u\) to \(C_j\times\{t\}\). Fix \(P\in D\). The integer \(\deg u_t^*\mathcal O_D(P)\) is independent of \(t\), because the curves \(C_j\times\{t\}\) are algebraically equivalent in \(C_j\times T\). It is positive for at least one \(t\), and hence for every \(t\). Thus every \(u_t\) is nonconstant.

Let \(d\) be this common positive degree. Riemann--Hurwitz applied to any \(u_t\) gives \(2g(C_j)-2\geq d(2g(D)-2)>0\), so \(g(C_j)\geq2\). By the universal property of the Hom scheme, \(u\) induces a morphism \(T\to\operatorname{Hom}_d(C_j,D)\), whose value at \(t\) is \(u_t\). The de Franchis theorem says that the set of nonconstant morphisms \(C_j\to D\) is finite \cite[Theorem~1, pp.~433--435]{HowardSommese}. Since \(T\) is reduced, its morphism to the Hom scheme factors through the reduced finite parameter space, and since \(T\) is connected, this morphism is constant. Therefore there is one nonconstant morphism \(v:C_j\to D\) such that \(u_t=v\) for every \(t\), which gives \(u=v\circ\operatorname{pr}_j\).
\end{proof}

\begin{lemma}\label{lem:product-lines-to-line}
Every nonconstant morphism \(u:(\mathbb P^1)^r\to\mathbb P^1\) factors through exactly one coordinate projection.
\end{lemma}

\begin{proof}
Put \(L=u^*\mathcal O_{\mathbb P^1}(1)\). There are uniquely determined integers \(d_i\geq0\) such that \(L\simeq\mathcal O(d_1,\ldots,d_r)\); here \(d_i\) is the degree of the restriction of \(u\) to a line obtained by varying only the \(i\)-th coordinate.

Suppose that \(d_a,d_b>0\) for distinct indices \(a,b\). Fix all other coordinates and denote the resulting map \(\mathbb P^1\times\mathbb P^1\to\mathbb P^1\) by \(u_{a,b}\). If \(P,Q\) are distinct points of the target, the divisors \(u_{a,b}^{-1}(P)\) and \(u_{a,b}^{-1}(Q)\) are disjoint. Both have numerical class \(d_aH_a+d_bH_b\), where \(H_a^2=H_b^2=0\) and \(H_aH_b=1\). Their intersection number must therefore be both \(0\) and \((d_aH_a+d_bH_b)^2=2d_ad_b>0\), a contradiction. Hence at most one \(d_i\) is positive. At least one is positive because \(u\) is nonconstant; denote its index by \(j\).

On every fiber of \(\operatorname{pr}_j:(\mathbb P^1)^r\to\mathbb P^1\), the restriction of \(u^*\mathcal O_{\mathbb P^1}(1)\) has degree zero on every coordinate line and is therefore trivial. Such a fiber is connected and projective, so its global regular functions are constant. The two sections of the trivial line bundle defining the restriction of \(u\) are consequently constant, and the restriction itself is constant. Fixing the other coordinates defines a morphism \(\phi:\mathbb P^1\to\mathbb P^1\), and constancy on the fibers gives \(u=\phi\circ\operatorname{pr}_j\). Finally, the unique positive entry of the multidegree is \(d_j=\deg\phi\), so no other coordinate can yield such a factorization.
\end{proof}

\begin{lemma}\label{lem:rational-partial-dependency}
Let \(Z\) be a possibly empty product of smooth projective connected curves of positive genus, and let \(u:(\mathbb P^1)^r\times Z\to\mathbb P^1\) be a morphism. Exactly one of the following occurs:
\begin{enumerate}[label=\textup{(\alph*)}]
\item \(u\) factors through \(Z\);
\item there is a unique index \(j\) for which \(u=\phi\circ\operatorname{pr}_j\), where \(\phi:\mathbb P^1\to\mathbb P^1\) has positive degree.
\end{enumerate}
\end{lemma}

\begin{proof}
For \(z\in Z\), write \(u_z\) for the restriction of \(u\) to \((\mathbb P^1)^r\times\{z\}\). The multidegree \((d_1,\ldots,d_r)\) of \(u_z\) is independent of \(z\): each \(d_i\) is the degree of the restriction of the fixed line bundle \(u^*\mathcal O_{\mathbb P^1}(1)\) to the \(i\)-th coordinate line. As \(z\) varies, these lines form the proper flat family \(\mathbb P^1\times Z\to Z\). The Hilbert polynomial of the restricted line bundle is locally constant in this family, and its linear coefficient is \(d_i\). Since \(Z\) is connected, every \(d_i\) is constant.

If every \(d_i\) is zero, then \(u_z\) is constant for every \(z\). After fixing \(x^0\in(\mathbb P^1)^r\), define \(v(z)=u(x^0,z)\). Fiberwise constancy gives \(u=v\circ\operatorname{pr}_Z\), so alternative \textup{(a)} holds.

Suppose that the multidegree is nonzero. By Lemma \ref{lem:product-lines-to-line}, it has a unique positive entry, say \(d_j=d\). The index \(j\) and the integer \(d\geq1\) are independent of \(z\). Fix the rational coordinates other than the \(j\)-th one. The restriction of \(u\) to the resulting copy of \(\mathbb P^1\times Z\) represents, by the universal property of the Hom scheme, a morphism \(Z\to\operatorname{Hom}_d(\mathbb P^1,\mathbb P^1)=\operatorname{Rat}_d\). For every \(z\), Lemma \ref{lem:product-lines-to-line} shows that \(u_z\) is the map represented by this restriction, composed with \(\operatorname{pr}_j\). The parameter space \(\operatorname{Rat}_d\) is affine \cite[Section~4.1]{SilvermanRat}. Since \(Z\) is connected and projective, every morphism from \(Z\) to an affine variety is constant. Thus the family is a single map \(\phi\), and \(u(x_1,\ldots,x_r,z)=\phi(x_j)\). The uniqueness of \(j\) follows from the multidegree, and the two alternatives are disjoint because \(d>0\) in \textup{(b)}.
\end{proof}

\subsection{Positive-genus factors}

\begin{lemma}\label{lem:positive-genus-stable-cover}
Let \(Z\) be a possibly empty product of smooth projective connected curves of positive genus, let \(T:Z\to Z\) be an endomorphism, and let \(q_0\geq1\). There are a multiple \(q\) of \(q_0\), an integer \(d\geq0\), a smooth projective variety \(Y\), a surjective \'{e}tale endomorphism \(\eta:Y\to Y\), and a surjective morphism \(\lambda_Z:Y\to Z_0:=(T^q)^d(Z)\) such that \(T^q(Z_0)=Z_0\) and \(T^q\circ\lambda_Z=\lambda_Z\circ\eta\).
\end{lemma}

\begin{proof}
Choose origins on the elliptic factors and write \(Z=H\times A\), where \(H=\prod_{i=1}^s C_i\) with \(g(C_i)\geq2\) for every \(i\), and \(A\) is a product of elliptic curves.

\noindent\emph{Step 1: the persistent hyperbolic coordinates.}
For each component of the \(H\)-coordinate of \(T\), Lemma \ref{lem:product-to-hyperbolic} gives either a constant map or a map that factors through one factor \(C_j\). A nonconstant factorization through an elliptic curve is impossible by Riemann--Hurwitz. We consequently obtain a partially defined function \(\tau\) on \(\{1,\ldots,s\}\): if the \(i\)-th output is nonconstant, then it depends only on the input \(C_{\tau(i)}\).

Let \(J\) be the set of indices on cycles of \(\tau\). The composite of the curve maps around any such cycle is a nonconstant self-map of a curve of genus at least \(2\). Riemann--Hurwitz shows that this composite has degree \(1\) and is an automorphism. Since its degree is the product of the degrees of the maps around the cycle, each of those maps is an isomorphism.

Choose \(q\) divisible by \(q_0\) and by every cycle length, and set \(G=T^q\). The preceding coordinate analysis shows that the \(H\)-coordinate of \(G(h,a)\) is a morphism \(\beta(h)\) independent of \(a\). Let \(G_A:H\times A\to A\) be the other coordinate and put \(v(h)=G_A(h,0)\). The morphism \(R(h,a)=G_A(h,a)-v(h)\) satisfies \(R(h,0)=0\). The rigidity theorem \cite[Theorem~2.1 and Corollaries~2.2 and~2.5, pp.~105--106]{MilneAbelian} implies first that \(R\) factors through \(A\) and then that the resulting origin-preserving morphism \(U:A\to A\) is a homomorphism. Consequently, \(G(h,a)=(\beta(h),U(a)+v(h))\). Let \(H_*=\prod_{j\in J}C_j\) and let \(p_H:H\to H_*\) be the projection. On every cyclic coordinate, \(\beta\) induces an automorphism; their product is an automorphism \(\gamma:H_*\to H_*\) satisfying \(p_H\circ\beta=\gamma\circ p_H\).

Apply Lemma \ref{lem:product-to-hyperbolic} to the coordinates of \(\beta\). Write the \(i\)-th coordinate as a constant \(c_i\) or as \(\rho_i\circ\operatorname{pr}_{\upsilon(i)}\), thereby defining a partial function \(\upsilon\). Since \(\beta\) is the \(H\)-coordinate of \(T^q\), the indices on cycles of \(\upsilon\) are exactly those in \(J\); the choice of \(q\) makes each of them a fixed point of \(\upsilon\). If \(i_0=i\) and \(i_{a+1}=\upsilon(i_a)\) as long as the latter is defined, direct induction gives two formulas. If \(i_0,\ldots,i_{n-1}\) all have successors, the \(i\)-th coordinate of \(\beta^n(h)\) is \(\rho_{i_0}\circ\cdots\circ\rho_{i_{n-1}}(h_{i_n})\). If \(i_a\) is the first index without a successor and \(a<n\), that coordinate is the constant \(\rho_{i_0}\circ\cdots\circ\rho_{i_{a-1}}(c_{i_a})\), with an empty composition interpreted as the identity.

Every path in this finite partial functional graph either enters a cycle or terminates. Define its transient depth to be the number of edges needed to reach \(J\), and to be one more than the number of edges preceding a terminal index in the terminating case. Let \(d_0\) be the maximum of these depths, with \(d_0=0\) if there are no indices. The formulas above show that, for \(n\geq d_0\), each coordinate of \(\beta^n(h)\) is either constant or determined by a coordinate \(h_j\) with \(j\in J\). Therefore \(p_H(h)=p_H(h')\) implies \(\beta^n(h)=\beta^n(h')\) for every \(n\geq d_0\). Fix a section \(s:H_*\to H\) of \(p_H\) by choosing a point on every factor not in \(J\).

\noindent\emph{Step 2: the stable abelian image and the stable image of \(G\).}
The images \(U^n(A)\) form a descending sequence of abelian subvarieties. Their dimensions eventually stabilize, and equality of the dimensions of two consecutive members forces equality. Choose \(m\geq0\) such that \(A_*:=U^m(A)=U^{m+1}(A)\). Then \(U_*:=U|_{A_*}\) is a surjective endomorphism of an abelian variety. It is an isogeny and, in characteristic zero, it is \'{e}tale.

For \(n\geq1\), define \(v_n(h)=\sum_{k=0}^{n-1}U^{n-1-k}v(\beta^k(h))\). Induction on \(n\), using \(v_{n+1}=U\circ v_n+v\circ\beta^n\), gives \(G^n(h,a)=(\beta^n(h),U^n(a)+v_n(h))\).

Each \(G^n(Z)\) is irreducible. It is closed because \(Z\) is projective and \(G^n:Z\to Z\) is proper, and \(G^{n+1}(Z)\subseteq G^n(Z)\). Once the dimensions of these images stabilize, two consecutive images are irreducible closed subvarieties of the same dimension and hence are equal. We may therefore choose \(d\geq\max\{1,d_0+m\}\) such that \(G^{d+1}(Z)=G^d(Z)\).

Define \(\lambda_H=\beta^d\circ s\). The choice of \(d_0\) gives \(\beta^d=\lambda_H\circ p_H\), and hence \(p_H\circ\lambda_H=\gamma^d\). Moreover, \(p_H(\beta(s(y)))=\gamma(y)=p_H(s(\gamma(y)))\). Applying \(\beta^d\) to these two points yields \(\beta\circ\lambda_H=\lambda_H\circ\gamma\). If \(\lambda_H(y)=\lambda_H(y')\), then applying \(p_H\) gives \(\gamma^d(y)=\gamma^d(y')\); thus \(y=y'\) because \(\gamma\) is an automorphism. In particular, \(\lambda_H\) is injective.

Put \(c=v_d\circ s\). For \(h\in H\), set \(y=p_H(h)\). If \(k\geq d_0\), then \(\beta^k(h)=\beta^k(s(y))\). If \(k<d_0\), then \(d-1-k\geq d-d_0\geq m\), so \(U^{d-1-k}(A)\subseteq U^m(A)=A_*\). Hence the terms with \(k\geq d_0\) cancel and
\[
 v_d(h)-c(y)=\sum_{k=0}^{d_0-1}U^{d-1-k}
 \bigl(v(\beta^k(h))-v(\beta^k(s(y)))\bigr)\in A_*.
\]
Give the closed image \(G^d(Z)\) its reduced induced structure and
denote it by \(Z_0\).
It is a closed subvariety of \(Z\), since \(G^d\) is proper.
The preceding calculation shows that
\[
 Z_0(\mathbb C)
 =
 \{(\lambda_H(y),c(y)+a):
 y\in H_*(\mathbb C),\ a\in A_*(\mathbb C)\}.
\]
Indeed, the preceding calculation proves the inclusion from left to right. For the reverse inclusion, given \((y,a)\in H_*\times A_*\), choose \(a_0\in A\) with \(U^d(a_0)=a\); then \(G^d(s(y),a_0)=(\lambda_H(y),c(y)+a)\). Our choice of \(d\) also gives \(G(Z_0)=Z_0\).

\noindent\emph{Step 3: an equivariant \'{e}tale cover.}
Define \(\lambda_Z:H_*\times A_*\to Z\) by \(\lambda_Z(y,a)=(\lambda_H(y),c(y)+a)\).
Since \(H_*\times A_*\) is reduced, the preceding pointwise
containment shows that
\(\lambda_Z:H_*\times A_*\to Z\) factors through \(Z_0\).
Its image contains every closed point of \(Z_0\), so the resulting
morphism \(H_*\times A_*\to Z_0\) is surjective.

For \(y\in H_*\), put \(\delta(y)=U(c(y))+v(\lambda_H(y))-c(\gamma(y))\). The point \(G(\lambda_Z(y,0))\) belongs to \(Z_0\), and its first coordinate is \(\beta(\lambda_H(y))=\lambda_H(\gamma(y))\). By the description of \(Z_0\), write \(G(\lambda_Z(y,0))=(\lambda_H(y'),c(y')+a_*)\) for some \(y'\in H_*\) and \(a_*\in A_*\). Equality of first coordinates and injectivity of \(\lambda_H\) give \(y'=\gamma(y)\); equality of second coordinates then gives \(\delta(y)=a_*\in A_*\). Since \(H_*\) is reduced and \(A_*\) is closed in \(A\), this pointwise containment implies that \(\delta:H_*\to A\) factors through \(A_*\).

Let \(Y=H_*\times A_*\) and define \(\eta(y,a)=(\gamma(y),U_*(a)+\delta(y))\). Direct substitution gives
\[
 \lambda_Z(\eta(y,a))
 =\bigl(\beta(\lambda_H(y)),U(c(y)+a)+v(\lambda_H(y))\bigr)
 =G(\lambda_Z(y,a)).
\]
Hence \(G\circ\lambda_Z=\lambda_Z\circ\eta\).

Finally, \(\eta\) is the composite of the finite \'{e}tale map \(\operatorname{id}_{H_*}\times U_*\), the automorphism \((y,a)\mapsto(y,a+\delta(y))\), and the automorphism \(\gamma\times\operatorname{id}_{A_*}\). It is therefore surjective and \'{e}tale. This proves every assertion of the lemma.
\end{proof}

\subsection{Proof of Mixed-Core Case: Theorem \ref{thm:mixed-core}}

\begin{proof}
Write \(B=R\times Z\), where \(R=(\mathbb P^1)^\ell\) and \(Z\) is the possibly empty product of the positive-genus factors.

\noindent\emph{Step 1: coordinate dependence.}
Every morphism \(\mathbb P^1\to C\) to a positive-genus curve is constant. Fix \(x^0\in R\) and define \(T(z)=\operatorname{pr}_Z(F(x^0,z))\). On every fiber of \(R\times Z\to Z\), vary the rational coordinates one at a time. Each component of \(\operatorname{pr}_Z\circ F\) is constant during every such variation, so it is constant on the whole fiber. Consequently, \(\operatorname{pr}_Z\circ F=T\circ\operatorname{pr}_Z\), and \(T:Z\to Z\) is an endomorphism. By Lemma \ref{lem:rational-partial-dependency}, each rational output coordinate either factors through \(Z\) or is \(\rho_i\circ\operatorname{pr}_{\sigma(i)}\) for a unique rational input coordinate \(\sigma(i)\) and a positive-degree map \(\rho_i:\mathbb P^1\to\mathbb P^1\). Thus \(\sigma\) is a partially defined function on the rational indices. Let \(J\) be the set of indices on cycles of \(\sigma\).

Let \(q_0\) be a common multiple of the cycle lengths, with \(q_0=1\) if \(J=\varnothing\). Apply Lemma \ref{lem:positive-genus-stable-cover} to \(T\) and \(q_0\). It gives a multiple \(q\) of \(q_0\), an integer \(d\), a stable image \(Z_0=(T^q)^d(Z)\), a smooth projective variety \(Y_0\), a surjective \'{e}tale endomorphism \(\eta_0:Y_0\to Y_0\), and a surjective morphism \(\lambda_Z:Y_0\to Z_0\) satisfying \(T^q\circ\lambda_Z=\lambda_Z\circ\eta_0\).

Set \(G=F^q\), let \(R_*=\prod_{j\in J}\mathbb P^1\), and let \(p_R:R\to R_*\) be the projection. Because \(q\) is divisible by every cycle length, the \(j\)-th coordinate of \(G\) is a self-map \(\theta_j:\mathbb P^1\to\mathbb P^1\) for every \(j\in J\). Let \(\theta=(\theta_j)_{j\in J}\). Each \(\theta_j\) has positive degree and is therefore surjective. With \(\operatorname{pr}_Z:B\to Z\) denoting the projection, we have \(p_R\circ G=\theta\circ p_R\) and \(\operatorname{pr}_Z\circ G=T^q\circ\operatorname{pr}_Z\).

\noindent\emph{Step 2: erasing the noncyclic rational inputs.}
If \(\sigma(i)\) is undefined, denote the corresponding morphism \(Z\to\mathbb P^1\) by \(\alpha_i\); otherwise the \(i\)-th rational coordinate of \(F\) is \(\rho_i\circ\operatorname{pr}_{\sigma(i)}\). Starting with \(i_0=i\), put \(i_{a+1}=\sigma(i_a)\) for as long as this is defined. Direct induction on \(n\geq1\) gives two forms for the \(i\)-th rational coordinate of \(F^n(x,z)\). If \(i_0,\ldots,i_{n-1}\) all have successors, it is \(\rho_{i_0}\circ\cdots\circ\rho_{i_{n-1}}(x_{i_n})\). If \(i_a\) is the first index without a successor and \(a<n\), it is \(\rho_{i_0}\circ\cdots\circ\rho_{i_{a-1}}(\alpha_{i_a}(T^{n-a-1}(z)))\), with an empty composition interpreted as the identity. This follows by applying the coordinate rule for \(F\) and using \(\operatorname{pr}_Z(F^k(x,z))=T^k(z)\).

Every path in the finite partial functional graph of \(\sigma\) either enters \(J\) or terminates, and it does so within at most \(\ell\) steps. Choose \(e\geq1\) with \(qe\geq\ell\). In the first formula with \(n=qe\), the final index \(i_n\) belongs to \(J\), so the value depends only on \(p_R(x)\). In the second formula it depends only on \(z\). The \(Z\)-coordinate of \(G^e=F^{qe}\) is \(T^{qe}(z)\) and also depends only on \(z\). We have therefore proved
\begin{equation}\label{eq:erasure}
 p_R(x)=p_R(x')
 \quad\Longrightarrow\quad
 G^e(x,z)=G^e(x',z)
 \quad\text{for all }x,x'\in R\text{ and }z\in Z.
\end{equation}

Fix a section \(s_R:R_*\to R\) of \(p_R\), and define \(\widehat\lambda:R_*\times Y_0\to B\) by \(\widehat\lambda(x,y)=G^e(s_R(x),\lambda_Z(y))\). Let \(\Pi=(p_R,\operatorname{pr}_Z):B\to R_*\times Z\). The identities at the end of Step 1 give \(\Pi\circ G^d=(\theta^d\circ p_R,(T^q)^d\circ\operatorname{pr}_Z)\). The displayed identity shows that \(\Pi\circ G^d\) is the product of
\(\theta^d\circ p_R\), which depends only on the \(R\)-coordinate,
and \((T^q)^d\circ\operatorname{pr}_Z\), which depends only on the
\(Z\)-coordinate.
Its image is therefore
\(\theta^d(R_*)\times(T^q)^d(Z)=R_*\times Z_0\), since
\(\theta\) is surjective. 

We claim that \(\widehat\lambda(R_*\times Y_0)=G^{d+e}(B)\). If \(w\in G^d(B)\), choose \(y\in Y_0\) with \(\lambda_Z(y)=\operatorname{pr}_Z(w)\). Equation \eqref{eq:erasure} gives \(G^e(w)=\widehat\lambda(p_R(w),y)\), proving \(G^{d+e}(B)\subseteq\widehat\lambda(R_*\times Y_0)\). Conversely, given \((x,y)\in R_*\times Y_0\), the equality \(\Pi(G^d(B))=R_*\times Z_0\) supplies \(w\in G^d(B)\) with \(p_R(w)=x\) and \(\operatorname{pr}_Z(w)=\lambda_Z(y)\). A second application of \eqref{eq:erasure} gives \(\widehat\lambda(x,y)=G^e(w)\in G^{d+e}(B)\). The claimed equality follows.

\noindent\emph{Step 3: equivariance and separation of the degree-one factors.}
Define \(\widehat\eta:R_*\times Y_0\to R_*\times Y_0\) by \(\widehat\eta(x,y)=(\theta(x),\eta_0(y))\). Put \(b_1=G(s_R(x),\lambda_Z(y))\) and \(b_2=(s_R(\theta(x)),\lambda_Z(\eta_0(y)))\). They have the same image under \(p_R\), because \(p_R\circ G=\theta\circ p_R\), and the same image under \(\operatorname{pr}_Z\), because \(T^q\circ\lambda_Z=\lambda_Z\circ\eta_0\). Equation \eqref{eq:erasure} therefore gives \(G^e(b_1)=G^e(b_2)\). By definition, this is \(G(\widehat\lambda(x,y))=\widehat\lambda(\widehat\eta(x,y))\), so \(G\circ\widehat\lambda=\widehat\lambda\circ\widehat\eta\).

Set \(B_0=G^{d+e}(B)=\widehat\lambda(R_*\times Y_0)\). Since every \(\theta_j\) and \(\eta_0\) is surjective, \(\widehat\eta\) is surjective, and equivariance gives \(G(B_0)=B_0\).

Partition \(J\) as \(J_{>1}\sqcup J_1\), according as \(\deg(\theta_j)\geq2\) or \(\deg(\theta_j)=1\). Put \(Y=(\mathbb P^1)^{J_1}\times Y_0\) and \(\eta=((\theta_j)_{j\in J_1},\eta_0)\). Each degree-one \(\theta_j\) is an automorphism of \(\mathbb P^1\), so \(\eta\) is a surjective \'{e}tale endomorphism of the smooth projective variety \(Y\). After reordering the factors, take \(\phi_1,\ldots,\phi_r\) to be the maps \(\theta_j\) with \(j\in J_{>1}\), and take \(\lambda=\widehat\lambda\) under the identification \(R_*\times Y_0\simeq(\mathbb P^1)^r\times Y\). Then every \(\deg(\phi_i)\geq2\), the required equivariance identity is the one already proved, and \(B_0=G^{d+e}(B)=(F^q)^{d+e}(B)\). Thus the theorem holds with \(N=d+e\).
\end{proof}

\section{Relatively \'{E}tale Systems over Products of Curves over \texorpdfstring{$\mathbb C$}{C}}
\label{sec:relative-curves-complex}

The stable model has a base of the form \((\mathbb P^1)^r\times Y\), with nonlinear split dynamics on the first factor and \'{e}tale dynamics on the second. The following observation allows the \'{e}tale factor to be absorbed into the relative directions, so that Theorem \ref{thm:relative-complex} applies without repeating the specialization argument.

\begin{lemma}[Forgetting an \'{e}tale base factor]\label{lem:forget-etale-factor}
Let \(R=(\mathbb P^1_{\mathbb C})^r\), let \(Y\) be a smooth quasi-projective complex variety, and put \(B=R\times Y\). Let \(\phi:R\to R\) be a split endomorphism whose coordinate maps have degree at least \(2\), let \(\eta:Y\to Y\) be \'{e}tale, and put \(\beta=(\phi,\eta)\). Suppose that \(p:X\to B\) is smooth, \(F:X\to X\) satisfies \(p\circ F=\beta\circ p\), and the relative graph of \(F\) over \(\beta\) is \'{e}tale. If \(q=\operatorname{pr}_R\circ p:X\to R\), then \(q\) is smooth and the relative graph of \(F\) over \(\phi\) is \'{e}tale.
\end{lemma}

\begin{proof}
The projection \(B=R\times Y\to R\) is smooth because \(Y\) is smooth, so its composite \(q\) with the smooth morphism \(p\) is smooth. Write \(p=(q,s)\), where \(s:X\to Y\), and put
\[
 Z_B=B\times_{\beta,B,p}X,
 \qquad
 Z_R=R\times_{\phi,R,q}X.
\]
The map \(\rho:Z_B\to Z_R\) that forgets the \(Y\)-coordinate is the base change of \(\eta:Y\to Y\) along \(Z_R\to Y\), \((r,x')\mapsto s(x')\). Indeed, \(Z_B\) is canonically \(Z_R\times_{s\circ\operatorname{pr}_2,Y,\eta}Y\). Hence \(\rho\) is \'{e}tale. If \(\delta_F^B\) and \(\delta_F^R\) denote the two relative graph morphisms, then \(\delta_F^R=\rho\circ\delta_F^B\). Both factors on the right are \'{e}tale, proving the assertion.
\end{proof}

\begin{theorem}[Relatively \'{e}tale systems over mixed bases]\label{thm:relative-mixed}
Let \(B=(\mathbb P^1_{\mathbb C})^r\times Y\) and \(\beta=(\varphi_1,\ldots,\varphi_r,\eta)\), where \(Y\) is smooth and quasi-projective, every \(\varphi_i:\mathbb P^1\to\mathbb P^1\) has degree at least \(2\), and \(\eta:Y\to Y\) is \'{e}tale. Let \(X\) be quasi-projective, let \(p:X\to B\) be smooth, and let \(F:X\to X\) satisfy \(p\circ F=\beta\circ p\). If the relative graph of \(F\) over \(\beta\) is \'{e}tale, then the dynamical Mordell--Lang conjecture holds for \(F\).
\end{theorem}

\begin{proof}
Apply Lemma \ref{lem:forget-etale-factor}. It makes \(F\) relatively \'{e}tale over \(\phi=(\varphi_1,\ldots,\varphi_r)\) through the smooth morphism \(q=\operatorname{pr}_{(\mathbb P^1)^r}\circ p\). Theorem \ref{thm:relative-complex} applies to this system and proves the result. When \(r=0\), the base in Theorem \ref{thm:relative-complex} is a point, exactly as stipulated there.
\end{proof}

We can now prove the relative theorem announced in the Introduction.

\begin{proof}[Proof of Theorem \ref{thm:relative-curve-base}]
Apply Theorem \ref{thm:mixed-core} to \(\beta\). Thus, for some \(q\geq1\) and \(N\geq0\), the closed subvariety \(B_0=(\beta^q)^N(B)\) satisfies \(\beta^q(B_0)=B_0\), and there are \(S=(\mathbb P^1)^r\times Y\), \(\psi=(\varphi_1,\ldots,\varphi_r,\eta)\), and a surjective morphism \(\lambda:S\to B_0\) satisfying \(\beta^q\circ\lambda=\lambda\circ\psi\). Put \(X_0=p^{-1}(B_0)\) and let \(p_0:X_0\to B_0\) be the restriction of \(p\). Since \(B_0\) is \(\beta^q\)-invariant, \(X_0\) is \(F^q\)-invariant. The variety \(X_0\) is closed in \(X\), hence quasi-projective, and \(p_0\) is smooth. Lemma \ref{lem:relative-etale-formal} first shows that \(F^q\) is relatively \'{e}tale over \(\beta^q\), and then, by equivariant base change along \(B_0\hookrightarrow B\), that \(F^q|_{X_0}\) is relatively \'{e}tale over \(\beta^q|_{B_0}\).

Form \(\widetilde X=S\times_{B_0}X_0\) and define \(\widetilde F(s,\xi)=(\psi(s),F^q(\xi))\). This is a self-map because \(\lambda(\psi(s))=\beta^q(\lambda(s))=\beta^q(p_0(\xi))=p_0(F^q(\xi))\). The projection \(\widetilde p:\widetilde X\to S\) is the base change of the smooth morphism \(p_0\), and hence is smooth. The variety \(\widetilde X\) is closed in \(S\times X_0\), because \(B_0\) is separated, and is therefore quasi-projective. The base-change assertion of Lemma \ref{lem:relative-etale-formal} shows that \(\widetilde F\) is relatively \'{e}tale over \(\psi\). Theorem \ref{thm:relative-mixed} consequently proves the dynamical Mordell--Lang conjecture for \(\widetilde F\).

Let \(\rho:\widetilde X\to X_0\) be the second projection. Given \(\xi\in X_0(\mathbb C)\), surjectivity of \(\lambda\) supplies \(s\in S(\mathbb C)\) with \(\lambda(s)=p_0(\xi)\); hence \((s,\xi)\in\widetilde X(\mathbb C)\) and \(\rho(s,\xi)=\xi\). Thus \(\rho\) is surjective on complex points and satisfies \(\rho\circ\widetilde F=(F^q|_{X_0})\circ\rho\). For \(\xi\in X_0(\mathbb C)\), choose \(\widetilde\xi\in\rho^{-1}(\xi)\). For every closed \(W\subseteq X_0\), the inverse image \(\rho^{-1}(W)\) is closed and \(S_{F^q|_{X_0}}(\xi,W)=S_{\widetilde F}(\widetilde\xi,\rho^{-1}(W))\). Therefore \(F^q|_{X_0}\) satisfies the conjecture.

It remains to pass from \(X_0\) to \(X\). Fix \(x\in X(\mathbb C)\) and a closed subvariety \(V\subseteq X\), and put \(x_0=(F^q)^N(x)\). Since \(p(x_0)=(\beta^q)^N(p(x))\in B_0\), we have \(x_0\in X_0\). With \(V_0=V\cap X_0\), the return set for \(F^q\) is
\[
 S_{F^q}(x,V)=
 \bigl(S_{F^q}(x,V)\cap\{n\in\mathbb N_0:n<N\}\bigr)
 \cup\bigl(N+S_{F^q|_{X_0}}(x_0,V_0)\bigr).
\]
The first set is finite, and the second is a finite union of arithmetic progressions. Hence the conjecture holds for \(F^q\) on \(X\), and Lemma \ref{lem:iteration} proves it for \(F\).
\end{proof}

\subsection{Immediate Consequences over \texorpdfstring{$\mathbb C$}{C}}
\label{sec:consequences}

We finish by proving the two consequences stated in the Introduction. In each case, the relative graph is identified explicitly before Theorem \ref{thm:relative-curve-base} is applied.

\begin{proof}[Proof of Theorem \ref{thm:all-curve-products}]
Take \(X=B\), \(p=\operatorname{id}_B\), and \(F=\beta\). The identity map \(p\) is smooth and \(p\circ F=\beta\circ p\). A point of \(B\times_{\beta,B,\operatorname{id}_B}B\) is a pair \((b,b')\) satisfying \(b'=\beta(b)\). Consequently, the morphism \(\alpha:B\to B\times_{\beta,B,\operatorname{id}_B}B\), given by \(\alpha(b)=(b,\beta(b))\), is an isomorphism with inverse the first projection. The relative graph \(\delta_F\) is exactly \(\alpha\) and is therefore \'{e}tale. Theorem \ref{thm:relative-curve-base} gives the result.
\end{proof}

\begin{proof}[Proof of Theorem \ref{thm:skew}]
Let \(p:B\times Y\to B\) be the first projection. Since \(Y\) is smooth, \(p\) is smooth; since \(B\) is projective and \(Y\) is quasi-projective, \(B\times Y\) is quasi-projective. The identity \(p\circ F=\beta\circ p\) follows from the definition of \(F\).

Every point of \(B\times_{\beta,B,p}(B\times Y)\) has the form \((b,(\beta(b),y'))\). The rule \(\iota(b,(\beta(b),y'))=(b,y')\) therefore defines an isomorphism to \(B\times Y\), whose inverse sends \((b,y')\) to \((b,(\beta(b),y'))\). For every \((b,y)\), one has \(\iota(\delta_F(b,y))=(b,\psi(b,y))\). Thus, under \(\iota\), the relative graph is the map assumed in the statement to be \'{e}tale. Theorem \ref{thm:relative-curve-base} applies and we have the theorem.
\end{proof}

\begin{proof}[Proof of Corollary \ref{cor:symmetric-powers}]
Put \(u=\operatorname{lcm}(u_1,\ldots,u_t)\) and \(r_j=u/u_j\), and
define \(\widetilde X=\prod_{j=1}^tX_j^{n_j}\),
\(\widetilde B=\prod_{j=1}^tB_j^{n_j}\),
\(\widetilde p=\prod_{j=1}^tp_j^{\times n_j}\), and
\(\widetilde\Theta=\prod_{j=1}^t\Theta_j^{\times n_j}\).
For each \(j\), apply Lemma \ref{lem:relative-etale-formal} to
\(F=\Theta_j^{u_j}\), \(p=p_j\), and \(\beta=\beta_j\), with iterate
\(r_j\). It follows that
\((\Theta_j^{u_j})^{r_j}=\Theta_j^u\) is relatively \'etale over
\(\beta_j^{r_j}\). Hence, under the canonical identification of the
corresponding fiber products, the relative graph of
\(\widetilde\Theta^u\) over
\(\widetilde\beta=
\prod_{j=1}^t(\beta_j^{r_j})^{\times n_j}\) is the product of these
relative graphs and is therefore \'etale. Since \(\widetilde p\) is
smooth, \(\widetilde X\) is quasi-projective, and \(\widetilde B\) is
a product of smooth projective connected curves, Theorem \ref{thm:relative-curve-base} applies
to \((\widetilde X,\widetilde B,\widetilde p,
\widetilde\Theta^u,\widetilde\beta)\). Thus
\(\widetilde\Theta^u\), and hence \(\widetilde\Theta\), satisfies
dynamical Mordell--Lang.
Let \(\pi:\widetilde X\to
\prod_{j=1}^t\operatorname{Sym}^{n_j}(X_j)\) be the product of the
symmetric-group quotient maps. It is finite, surjective, and
equivariant for \(\widetilde\Theta\) and
\(\prod_{j=1}^t\operatorname{Sym}^{n_j}(\Theta_j)\). Dynamical
Mordell--Lang descends through such a morphism by lifting the initial
point and pulling back the target subvariety, proving the first
assertion. For the final assertion, take \(B_j=X_j\),
\(p_j=\operatorname{id}_{X_j}\), \(u_j=1\), and
\(\beta_j=\Theta_j\); the relative graph is then an isomorphism.
\end{proof}

\bibliographystyle{amsplain}
\bibliography{Rel_Etale_Final_Arxiv}

\end{document}